# CONVERGENCE RATES FOR POSTERIOR DISTRIBUTIONS AND ADAPTIVE ESTIMATION

By Tzee-Ming Huang

*Iowa State University*

The goal of this paper is to provide theorems on convergence rates of posterior distributions that can be applied to obtain good convergence rates in the context of density estimation as well as regression. We show how to choose priors so that the posterior distributions converge at the optimal rate without prior knowledge of the degree of smoothness of the density function or the regression function to be estimated.

**1. Introduction.** Bayesian methods have been used for nonparametric inference problems, and many theoretical results have been developed to investigate the asymptotic properties of nonparametric Bayesian methods. So far, the positive results are on consistency and convergence rates. For example, Doob (1949) proved the consistency of posterior distributions with respect to the joint distribution of the data and the prior under some weak conditions, and Schwartz (1965) extended Doob's result to Bayes decision procedures with possibly nonconvex loss functions. For the frequentist version of consistency, see Diaconis and Freedman (1986) for a review on consistency results on tail-free and Dirichlet priors. Barron, Schervish and Wasserman (1999) gave some conditions to achieve the frequentist version of consistency in general. Ghosal, Ghosh and Ramamoorthi (1999) also gave a similar consistency result and applied it to Dirichlet mixtures.

For convergence rates, there are some general results by Ghosal, Ghosh and van der Vaart (2000) and Shen and Wasserman (2001). However, there are few results on adaptive estimation in the study of posterior convergence rates. Belitser and Ghosal (2003) dealt with adaptive estimation in the infinite normal mean set-up. In this paper, we also have results on adaptive estimation, but these are done in the density estimation and regression setups.

The goal of this paper is to develop theorems on convergence rates for posterior distributions which can be used for adaptive estimation. In this paper we have theorems on convergence rates in two contexts: density estimation and regression. In either case, we consider the Bayesian estimation of some function $f$ (a density function or a regression function) based on a sample $(Z_1, \ldots, Z_n)$ and are interested in the convergence rates for the posterior distributions for $f$.

Below is the specific problem setup. Suppose that when $f$ is given, $(Z_1, \ldots, Z_n)$ is a random sample from a distribution with density $p_f$ with respect to a measure $\mu$ on a sample space $(\mathcal{S}, \mathcal{B})$, $f_o$ is the true value for $f$, and $f_o$ belongs to some function space $\mathcal{F}$. Suppose that $\tilde{\pi}$ is a prior on $\mathcal{F}$ and $\tilde{B}_d(s_n) = \{f \in \mathcal{F} : d(f, f_o) \leq s_n\}$ is an $s_n$ neighborhood of $f_o$ with respect to the metric $d$, where $d$ is the Hellinger distance in the density estimation case and is the $L_2$ distance in the regression case.

We would like to show that the posterior probability

$$(1) \qquad \tilde{\pi}(\tilde{B}_d(s_n)^c | Z_1, \ldots, Z_n) = \frac{\int_{\tilde{B}_d(s_n)^c} \prod_{i=1}^n p_f(Z_i) \, d\tilde{\pi}(f)}{\int_{\mathcal{F}} \prod_{i=1}^n p_f(Z_i) \, d\tilde{\pi}(f)}$$

converges to zero in $P_{f_o}^n$ probability, and the rate $s_n$ is as good as if the degree of smoothness of $f_o$ were known. This is known as the adaptive estimation problem.

For the purpose of adaptive estimation, we take $\mathcal{F}$ to be $\bigcup_{j \in J} \mathcal{F}_j$, where $J$ is a countable index set (not necessarily a set of integers) and the $\mathcal{F}_j$'s are function spaces of different degrees of smoothness. A natural way to construct priors on $\mathcal{F}$ is to consider sieve priors. A sieve prior is a prior $\tilde{\pi}$ of the following form:

$$\tilde{\pi} = \sum_{j \in J} a_j \tilde{\pi}_j,$$

where $a_j \geq 0$, $\sum_{j \in J} a_j = 1$, and each $\tilde{\pi}_j$ is a prior defined on $\mathcal{F}$ but supported on $\mathcal{F}_j$. To make it easier to specify the $\tilde{\pi}_j$'s, we assume that each $\mathcal{F}_j$ is finite-dimensional and can be represented as $\{f_{\theta,j} : \theta \in \Theta_j\}$ for some parameter space $\Theta_j$. We also assume that each $\tilde{\pi}_j$ is induced by a prior $\pi_j$ defined on $\Theta_j$. Then the posterior probability in (1) can be written as $U_n/V_n$, where

$$U_n = \sum_j a_j \int_{B_{d,j}(s_n)^c} \prod_{i=1}^n \frac{p_{f_{\theta,j}}(Z_i)}{p_{f_o}(Z_i)} \, d\pi_j(\theta)$$

and

$$V_n = \sum_j a_j \int_{\Theta_j} \prod_{i=1}^n \frac{p_{f_{\theta,j}}(Z_i)}{p_{f_o}(Z_i)} \, d\pi_j(\theta),$$

where $B_{d,j}(s_n) = \{\theta \in \Theta_j : d(f_{\theta,j}, f_o) \leq s_n\}$.

This paper is organized as follows. Section 2 gives a theorem on convergence rates in the density estimation case and some examples of applying the theorem to obtain adaptive rates. Section 3 contains the same things as in Section 2, but in the context of regression. Proofs are in Section 4.

## 2. Density estimation.

2.1. *Theorem.* This section gives a convergence rate theorem for Bayesian density estimation. The setup is as described in Section 1, with $p_f = f$ and $d$ being the Hellinger metric $d_H$, which is defined by

$$d_H(f, g) = \sqrt{\int (\sqrt{f} - \sqrt{g})^2 \, d\mu}.$$

To make the posterior probability $U_n/V_n \to 0$, we need some conditions to give bounds for $U_n$ and $V_n$.

To bound $U_n$, we will make an assumption about the structure of each parameter space $\Theta_j$, and then specify the $a_j$ accordingly. Let $\|\cdot\|_\infty$ denote the sup norm

$$B_{d_H,j}(\eta, r) = \{\theta \in \Theta_j : d_H(f_{\eta,j}, f_{\theta,j}) \leq r\}$$

and $N(B, \delta, d')$ denote the $\delta$-covering number of a set $B$ with respect to a metric $d'$, which is defined as the smallest number of $\delta$-balls (with respect to $d'$) that are needed to cover the set $B$. Here is the assumption.

ASSUMPTION 1. *For each $j \in J$, there exist constants $A_j$ and $m_j$ such that $A_j \geq 0.0056$, $m_j \geq 1$, and for any $r > 0$, $\delta \leq 0.0056r$, $\theta \in \Theta_j$,*

$$N(B_{d_H,j}(\theta, r), \delta, d_{j,\infty}) \leq \left(\frac{A_j r}{\delta}\right)^{m_j},$$

*where $d_{j,\infty}(\theta, \eta)$ is defined as $\|\log f_{\theta,j} - \log f_{\eta,j}\|_\infty$ for all $\theta, \eta \in \Theta_j$.*

Suppose Assumption 1 holds. We specify the $a_j$'s in the following way:

$$(2) \qquad a_j = \alpha \exp\left(-\left(1 + \frac{1-4\gamma}{8}\right)\eta_j\right),$$

where $\alpha$ is a normalizing constant so that $\sum_j a_j = 1$, $\gamma \doteq 0.1975$ is the solution to $0.13\gamma/\sqrt{1-4\gamma} = 0.0056$, and

$$(3) \qquad \eta_j = \frac{4m_j}{1-4\gamma} \log\left(\frac{46.2 A_j \sqrt{1-4\gamma}}{\gamma}\right) + \frac{8C_j}{1-4\gamma}$$

for some $C_j$ such that $C_j \geq 0$ and $\sum_j e^{-C_j} \leq 1$.

Note:

1. Assumption 1 is based on Assumption 1 in Yang and Barron (1998) so that their results can be applied here. The constants $A_j$ and $m_j$ can be figured out based on the local structure of $\Theta_j$. In many cases, $m_j$ can be taken as the dimension of $\Theta_j$, as stated in Lemma 2.
2. The constants $C_j$'s are here to make sure that $\sum_j a_j < \infty$ since $a_j \leq \alpha e^{-C_j}$. Indeed, we may take $\eta_j$ to be some large constant times $m_j \log A_j$, if this choice makes $\{a_j\}$ summable. Also, specific constant values are given in (2) and (3) for calculational convenience. Different choices are possible.

To find a bound for $V_n$, we will use Lemma 1 of Shen and Wasserman (2001), which says we can bound $V_n$ from below if the prior puts enough probability on a small neighborhood of the true density $f_o$. To guarantee enough prior probability around $f_o$, we proceed as follows.

1. Find a model $\mathcal{F}_{j_n}$ that receives enough weight $a_{j_n}$ and is close to $f_o$, that is, there exists $\beta_n$ in $\Theta_{j_n}$ so that $f_{\beta_n, j_n}$ is close to $f_o$.
2. Make sure the prior $\pi_{j_n}$ puts enough probability on a neighborhood of $\beta_n$. This helps $\tilde{\pi}$ put some probability around $f_o$ since $a_{j_n}$ is not too small.

For the first step, we simply assume that it is possible.

ASSUMPTION 2. There exist $j_n$ and $\beta_n \in \Theta_{j_n}$ such that

$$(4) \qquad \max(D(f_o\|f_{\beta_n,j_n}), V(f_o\|f_{\beta_n,j_n})) + \frac{\eta_{j_n}}{n} \leq \varepsilon_n^2$$

for some sequence $\{\varepsilon_n\}$, where $D(f\|g) = \int f \log(f/g)\, d\mu$, $V(f\|g) = \int f(\log(f/g))^2\, d\mu$, $\eta_{j_n}$ is as defined in (3) with $A_{j_n}$ and $m_{j_n}$ in Assumption 1.

Before going to assumptions for the second step, we add one more condition here to allow us to use neighborhoods that are different but comparable to the neighborhoods in Lemma 1 of Shen and Wasserman (2001).

ASSUMPTION 3. For the $j_n$ in Assumption 2, there exists a metric $d_{j_n}$ on $\Theta_{j_n}$ such that

$$(5) \qquad \int f_o \left(\log \frac{f_{\eta,j_n}}{f_{\theta,j_n}}\right)^2 d\mu \leq K_0' d_{j_n}^2(\eta, \theta)$$

for all $\eta$, $\theta$ in $\Theta_{j_n}$, and

$$D(f_o\|f_{\theta,j_n}) \leq K_0'' V(f_o\|f_{\theta,j_n})$$

for all $\theta \in \Theta_{j_n}$, where $K_0'$ and $K_0''$ are constants independent of $n$.

The following two assumptions are for the second step.

ASSUMPTION 4. *For $j_n$, $A_{j_n}$, $m_{j_n}$, $\beta_n$, $\varepsilon_n$ and $d_{j_n}$ in Assumptions 1–3, there exists $b_1 \geq 0$ such that*

$$N(\Theta_{j_n}, \varepsilon_n, d_{j_n}) \leq (A_{j_n}^{b_1} K_4)^{m_{j_n}},$$

*where $N(\Theta_{j_n}, \varepsilon_n, d_{j_n})$ is the $\varepsilon_n$-covering number of $\Theta_{j_n}$ with respect to the metric $d_{j_n}$.*

ASSUMPTION 5. *For $j_n$, $A_{j_n}$, $m_{j_n}$, $\beta_n$, $\varepsilon_n$ and $d_{j_n}$ in Assumptions 1–3, there exist constants $K_5$ and $b_2 \geq 0$ such that for any $\theta_1 \in \Theta_{j_n}$,*

$$\frac{\pi_{j_n}(B_{d_{j_n}, j_n}(\theta_1, \varepsilon_n))}{\pi_{j_n}(B_{d_{j_n}, j_n}(\beta_n, \varepsilon_n))} \leq (A_{j_n}^{b_2} K_5)^{m_{j_n}}.$$

Note:

1. Assumption 4 is here to give more control of the overall size of $\Theta_{j_n}$ in terms of the $\varepsilon_n$-covering number (Assumption 1 essentially deals with the local structure). This control is to prevent the total prior probability from getting spread out so much that each neighborhood gets little probability.
2. Assumption 5 is to make sure that the prior supported on $\Theta_{j_n}$ puts enough probability near $\beta_n$ compared to some other neighborhood.

Finally, we assume the following.

ASSUMPTION 6. *As $n \to \infty$,*

$$\varepsilon_n \to 0 \quad \text{and} \quad n\varepsilon_n^2 \to \infty.$$

Now we have the following theorem.

THEOREM 1. *Suppose that Assumptions 1–6 hold. Then with $a_j$ defined in (2), there exist positive constants $c$, $K_1$ and $K_2$ that are independent of $n$, so that*

(6) $$\tilde{\pi}(\tilde{B}_{d_{\mathrm{H}}}(K_1\varepsilon_n)^c | X_1, \ldots, X_n) \leq c\exp(-K_2 n\varepsilon_n^2)$$

*except on a set of probability converging to zero.*

The proof of Theorem 1 is given in Section 4.

2.2. *Example*: *spline basis*. In this section, we assume that $\log f_o$ is in the Sobolev space $W_\infty^s[0,1] = \{g : \|D^s g\|_{L_\infty[0,1]} < \infty\}$, where $s$ is a positive integer and $\|\cdot\|_{L_\infty[0,1]}$ is the essential sup norm with respect to the Lebesgue measure on $[0,1]$. We will see that using the sieve prior given below, the posterior distribution converges at the rate $n^{-s/(1+2s)}$ in Hellinger distance.

LEMMA 1. *Suppose that* $\log f_o \in W^s_\infty[0,1]$ *as defined above and $\mu$ is the Lebesgue measure on* $[0,1]$. *Let* $J = \{(k,q,L) : k, q$ *and $L$ are integers* $k \geq 0, q \geq 1,$ *and* $L \geq 1\}$. *For* $j = (k,q,L) \in J$, *let* $m_j = k + q$, *and for* $i \in \{1, \ldots, m_j\}$, *let* $B_{j,i}$ *be the normalized B-spline associated with the knots* $y_i, \ldots, y_{i+q}$ *as in Definition* 4.19, *page* 124 *in Schumaker ([1981](#)), where*

$$(y_1, \ldots, y_q, y_{q+1}, \ldots, y_{q+k}, y_{q+k+1}, \ldots, y_{2q+k})$$
$$= (\underbrace{0, \ldots, 0}_{q \text{ times}}, 1/(1+k), \ldots, k/(1+k), \underbrace{1, \ldots, 1}_{q \text{ times}}).$$

*Define*

$$\Theta_j = \{\theta \in R^{m_j} : \theta' \mathbb{1}_{m_j} = 0, \|D^r \log f_{\theta,j}\|_{L_\infty[0,1]} \leq L, \forall r \in \{0, 1, \ldots, q-1\}\},$$

*where* $\mathbb{1}_{m_j} = (1, \ldots, 1)' \in R^{m_j}$, $\log f_{\theta,j} = -\psi(\theta) + \theta' B$, $\psi(\theta) = \log \int_0^1 e^{\theta' B(x)} dx$ *is the normalizing constant, and* $B = (B_{j,1}, \ldots, B_{j,m_j})$. *Define* $\eta_j$ *as in* [(3)](#) *with*

(7) $\quad A_j = 19.28\sqrt{q}(2q+1)9^{q-1}(L+1)e^{L/2} + 0.06 \quad and \quad C_j = m_j + L;$

*define $a_j$ as in* [(2)](#). *Let $\pi_j$ be the Lebesgue measure on $\Theta_j$. Let $\tilde{\pi}_j$ be the induced prior of $\pi_j$ and $\tilde{B}_{d_H}(s_n)$ denote the $s_n$ Hellinger neighborhood of $f_o$, as defined on page* 3 *of Schumaker ([1981](#)). Then for the prior* $\tilde{\pi} = \sum_j a_j \tilde{\pi}_j$, *the posterior probability* $\tilde{\pi}(\tilde{B}_{d_H}(s_n)^c | X_1, \ldots, X_n)$ *converges to zero in probability for some* $s_n \propto n^{-s/(1+2s)}$.

The proof of Lemma [1](#) is given in Section [4](#).
Note:

1. Log-spline models have been used in density estimation and give good convergence rates; see Stone ([1990](#)), for example.
2. The prior does not depend on $s$, but it adapts to the smoothness parameter $s$.
3. Here we take $\pi_j$ to be the Lebesgue measure on $\Theta_j$, but we may also take $\pi_j$ to be some measure that has a density $q_j$ with respect to the Lebesgue measure on $\Theta_j$. As long as $\|\log q_j\|_\infty$ is uniformly bounded in $j$, the convergence rates should be the same.
4. $C_j = m_j + L$ is just one possible choice. In general, if we choose $\{C_j\}$ so that $\sum_j e^{-C_j} < \infty$ and $C_{j_n} \to \infty$ no faster than $m_{j_n} \log A_{j_n}$, where $j_n$ is as in Assumption [2](#), then it should be a good choice.
5. To figure out $A_j$ and $m_j$, the following lemma, from Lemma 1 by Yang and Barron ([1998](#)), is useful.

LEMMA 2. *Suppose that* $\{S_l : l \in \Lambda\}$ *is a countable collection of linear function spaces on* $[0,1]$. *Suppose that for each $S_l$ there is a basis* $\{B_{l,1}, \ldots, B_{l,m_l}\}$.

Suppose that there exist constants $T_1$ and $T_2$ such that for $\theta = (\theta_1, \ldots, \theta_{m_l}) \in R^{m_l}$,

(8) $$\left\|\sum_{i=1}^{m_l} \theta_i B_{l,i}\right\|_\infty \leq T_1 \max_i |\theta_i|$$

and

(9) $$\left\|\sum_{i=1}^{m_l} \theta_i B_{l,i}\right\|_2 \geq \frac{T_2}{\sqrt{m_l}} \sqrt{\sum_{i=1}^{m_l} \theta_i^2},$$

where $\|\cdot\|_2$ denotes the $L_2$ norm with respect to the Lebesgue measure on $[0,1]$. Let

(10) $$\log f_{\theta,j} = -\psi(\theta) + \sum_{i=1}^{m_l} \theta_i B_{l,i},$$

where $\psi(\theta) = \log \int_0^1 \exp(\sum_{i=1}^{m_l} \theta_i B_{l,i}(x)) \, dx$ is the normalizing constant. Suppose that $1 \in S_l$ for all $l \in \Lambda$, $J = \{(l, L) : l \in \Lambda, L \text{ is a positive integer}\}$ and for $j \in J$,

$$\Theta_j \subset \{\theta \in R^{m_l} : \|\log f_{\theta,j}\|_\infty \leq L\}.$$

Then Assumption 1 holds with

(11) $$A_j = 19.28 \frac{T_1}{T_2}(L+1)e^{L/2} + 0.06 \quad \text{and} \quad m_j = m_l.$$

2.3. *Example*: *Haar basis*. In this section, we assume that $\log f_o$ is a continuous function on $[0,1]$ with $\|\log f_o\|_\infty \leq M_0$, and we approximate $\log f_o$ using the Haar basis $\{\mathbb{1}_{[0,1]}(x), \psi_{j_1,k_1}(x) : 0 \leq j_1, 0 \leq k_1 \leq 2^{j_1} - 1\}$, where $\psi_{j_1,k_1}(x) = 2^{j_1/2}\psi^*(2^{j_1}x - k_1)$ and $\psi^*(x) = \mathbb{1}_{[0,0.5]}(x) - \mathbb{1}_{[0.5,1]}(x)$. We also assume that the coefficients of the $L_2$ expansion of $\log f_o$ for the Haar basis, denoted by $d_{j_1,k_1}$, satisfy the following condition:

(12) $$\sum_{j_1 \geq 0} (2^{j_1+1} - 1)^{2\alpha} \sum_{k_1=0}^{2^{j_1}-1} d_{j_1,k_1}^2 \leq H_0^2$$

for some $H_0 > 0$ and $\alpha \in (0,1)$. According to Barron, Birgé and Massart [(1999), page 330], the above condition on the Haar basis coefficients corresponds to the Besov space $B_{2,2}^\alpha[0,1]$. The Besov space $B_{2,2}^\alpha[0,1]$ is indeed the Sobolev space $W_2^\alpha[0,1]$, so the optimal convergence rate is $n^{-\alpha/(1+2\alpha)}$ in $L_2$-distance. We will see that using the sieve prior given below, the posterior distribution converges at the rate $n^{-\alpha/(1+2\alpha)}(\log n)^{1/2}$ in Hellinger distance, which is close to the optimal rate $n^{-\alpha/(1+2\alpha)}$ within a $(\log n)^{1/2}$ factor:

LEMMA 3. *Suppose that $\log f_o$ is in the space specified above and $\mu$ is the Lebesgue measure on $[0,1]$. Let $J = \{(l, L): l \text{ and } L \text{ are integers. } l \geq 0, L \geq 1\}$. For $j = (l, L) \in J$, let $m_j = 2^{l+1}$. Reindex the Haar basis in the following way:*

$$\{\psi_{j_1,k_1}: 0 \leq j_1 \leq l, 0 \leq k_1 \leq 2^{j_1} - 1\} \stackrel{def}{=} \{B_{j,i}: 1 \leq i \leq m_j - 1\}.$$

*Then for $\theta \in R^{m_j-1}$, define $\log f_{\theta,j} = -\psi(\theta) + \theta'B$, where $\psi(\theta) = \log \int_0^1 e^{\theta'B(x)}\,dx$ is the normalizing constant and $B = (B_{j,1}, \ldots, B_{j,m_j})$. Define*

$$\Theta_j = \{\theta \in R^{m_j-1}: \|\theta'B\|_\infty \leq L\}$$

*and let $\pi_j$ be the Lebesgue measure on $\Theta_j$. Define $a_j$ and $\eta_j$ according to (2) and (3) with*

(13) $\qquad A_j = 19.28 \cdot 2^{(l+1)/2}(2L+1)e^L + 0.06 \quad \text{and} \quad C_j = m_j + L.$

*Let $\pi_j$ be the Lebesgue measure on $\Theta_j$. Let $\tilde{\pi}_j$ be the induced prior of $\pi_j$ and $\tilde{B}_{d_H}(s_n)$ denote the $s_n$ Hellinger neighborhood of $f_o$, as defined on page 3 in Schumaker (1981). Then for the prior $\tilde{\pi} = \sum_j a_j \tilde{\pi}_j$, the posterior probability $\tilde{\pi}(\tilde{B}_{d_H}(s_n)^c | X_1, \ldots, X_n)$ converges to zero in probability for some $s_n \propto n^{-\alpha/(1+2\alpha)}(\log n)^{1/2}$.*

The proof of Lemma 3 is given in Section 4.
Note:

1. For the choice of $a_j$ and $\pi_j$, see the note for Lemma 1.
2. To specify $A_j$ and $m_j$, Lemma 2 is no longer applicable since $T_1$ in (8) cannot be taken as a constant in this case. We use the following lemma [from Lemma 2 by Yang and Barron (1998)] instead.

LEMMA 4. *Suppose that $\{S_l : l \in \Lambda\}$ is a countable collection of linear function spaces on $[0,1]$ and that for each $l$ there exists a constant $K_l > 0$ such that for all $h \in S_l$,*

(14) $$\|h\|_\infty \leq K_l \|h\|_2.$$

*Suppose that each $S_l$ is spanned by a bounded and linearly independent (under $L_2$ norm) basis $1, B_{l,1}, \ldots, B_{l,m_l}$. For $\theta \in R^{m_l}$, define $\log f_{\theta,j} = -\psi(\theta) + \sum_{i=1}^{m_l} \theta_i B_{l,i}$, where $\psi(\theta) = \log \int_0^1 \exp(\sum_{i=1}^{m_l} \theta_i B_{l,i}(x))\,dx$. Suppose that $J = \{(l, L): l \in \Lambda, L \text{ is a positive integer}\}$ and for each $j \in J$,*

(15) $$\Theta_j \subset \{\theta \in R^{m_l}: \|\log f_{\theta,j}\|_\infty \leq 2L\}.$$

*Then Assumption 1 holds with*

(16) $\qquad A_j = 19.28 K_l(2L+1)e^L + 0.06 \quad \text{and} \quad m_j = m_l + 1.$

In the spline density estimation result, the convergence rate is optimal and we have full adaption. But the Haar basis result here is quite different. The convergence rate involves an extra log factor, which comes from the $K_l$ in (16). In the spline case there is no $K_l$ and $A_j$ is approximately a constant when $j = j_n$ for large $n$ ($j_n$ is the index for one of the best models at sample size $n$). In this case $A_j$ is approximately proportional to the model dimension $m_j$ when $j = j_n$ because of the factor $K_l$.

## 3. Regression.

3.1. *Theorem.* In this section, a Bayesian convergence rate theorem is given in the context of regression. The setup is as described in Section 1, with $Z_i = (X_i, Y_i)$, where $Y_i = f(X_i) + \varepsilon_i$, $X_i$ and $\varepsilon_i$ are independent, $X_i$ is distributed according to some probability measure $\mu_X$ and $\varepsilon_i$ is normally distributed with mean zero and known variance $\sigma^2$. Thus the density $p_f$ (with respect to $\mu_X \times$ Lebesgue measure on $R$) is

$$p_f(x,y) = \frac{1}{\sqrt{2\pi}\sigma} e^{-(y-f(x))^2/(2\sigma^2)}.$$

The metric $d$ is the $L_2(\mu_X)$ metric. We also assume that $\|f_o\|_\infty$ is bounded by a known constant $M$.

To bound $U_n$ and $V_n$, we modify the assumptions in Theorem 1 in the following way. Let

$$B_{L_2(\mu_X),j}(\eta,r) = \{\theta \in \Theta_j : \|f_{\eta,j} - f_{\theta,j}\|_{L_2(\mu_X)} \leq r\}.$$

Assumption 1 is replaced with the following.

ASSUMPTION 7. For each $j$, there exist constants $A_j$ and $m_j$ such that $0 < A_j \leq 0.0056$, $m_j \geq 1$, and for any $r > 0$, $\delta \leq 0.0056r$, $\theta \in \Theta_j$,

$$N(B_{L_2(\mu_X),j}(\theta,r),\delta,d_{j,\infty}) \leq \left(\frac{A_j r}{\delta}\right)^{m_j},$$

where $d_{j,\infty}(\theta,\eta) = \|f_{\theta,j} - f_{\eta,j}\|_\infty$ for all $\theta, \eta \in \Theta_j$.

Also, suppose Assumption 7 holds: we specify the weights $a_j$ in the following way to give an upper bound for $U_n$:

(17) $$a_j = \alpha \exp\left(-\left(1 + \frac{1}{2\sigma^2} + \frac{0.0056}{\sigma}\right)\eta_j\right),$$

where $\alpha$ is a normalizing constant so that $\sum_j a_j = 1$ and

(18) $$\eta_j = \frac{4m_j}{c_{1,M,\sigma}(1-4\gamma)} \log(1072.5 A_j) + C_j \max\left(1, \frac{8}{c_{1,M,\sigma}(1-4\gamma)}\right)$$

for some $C_j$ such that $C_j \geq 0$ and $\sum_j e^{-C_j} \leq 1$.

Assumption 2 is replaced with the following assumption.

ASSUMPTION 8. *There exist $j_n$ and $\beta_n \in \Theta_{j_n}$ such that*

$$\max(D(p_{f_o}\|p_{f_{\beta_n,j_n}}), V(p_{f_o}\|p_{f_{\beta_n,j_n}})) + \frac{\eta_{j_n}}{n} \leq \varepsilon_n^2 \qquad (19)$$

*for some sequence $\{\varepsilon_n\}$, where $\eta_{j_n}$ is as defined in* (18) *with $A_{j_n}$ and $m_{j_n}$ in Assumption* 7.

Assumption 3 is replaced with the following.

ASSUMPTION 9.

$$\|f_{\theta,j_n} - f_{\eta,j_n}\|^2_{L_2(\mu_X)} \leq K_0' d_{j_n}^2(\theta, \eta) \qquad \text{for all } \theta, \eta \in \Theta_{j_n}. \qquad (20)$$

Assumptions 4–6 remain unchanged except that "Assumptions 1–3" should be changed to "Assumptions 7–9."

Now we have the following theorem.

THEOREM 2. *Suppose that $\|f_{\theta,j}\|_\infty \leq M$ for all $j$ and $\theta \in \Theta_j$. Suppose that Assumptions* 7–9 *and Assumptions* 4–6 *hold with the reference change made as mentioned above. Then with $a_j$ defined in* (17), *there exists a positive constant $K_1$ such that $\tilde{\pi}(\tilde{B}_{L_2(\mu_X)}(K_1\varepsilon_n)^c|X_1,\ldots,X_n)$ converges to zero in probability. Here $\tilde{B}_{L_2(\mu_X)}(K_1\varepsilon_n)$ denotes the $K_1\varepsilon_n$ neighborhood of $f_o$ with respect to the $L_2(\mu_X)$ metric, as defined on page* 1557.

The proof of Theorem 2 is given in Section 4.

3.2. *An example.* In this section, we consider $f_o \in W_\infty^s[0,1] = \{g : \|D^s g\|_{L_\infty[0,1]} < \infty\}$ and approximate $f_o$ using a spline basis. The minimax rate for this space in $L_2$ metric, according to Stone (1982), is $n^{-s/(1+2s)}$. We will see that, using the sieve prior given below, the posterior distribution converges at the optimal rate $n^{-s/(1+2s)}$ in $L_2$ distance.

LEMMA 5. *Suppose that $f_o \in W_\infty^s[0,1]$, $\|f_o\|_\infty < M$, where $M$ is a known constant. Suppose that $\mu_X$ is the Lebesgue measure on $[0,1]$. Let $J = \{(k,q,L) : k,q \text{ and } L \text{ are integers}; k \geq 0, q \geq 1, L \geq 1\}$. For $j = (k,q,L) \in J$, let $m_j = k + q$, and for $i \in \{1,\ldots,m_j\}$, let $B_{j,i}$ be the normalized B-spline associated with the knots $y_i,\ldots,y_{i+q}$, where*

$$(y_1,\ldots,y_q, y_{q+1},\ldots,y_{q+k}, y_{q+k+1},\ldots,y_{2q+k})$$
$$= (\underbrace{0,\ldots,0}_{q \text{ times}}, 1/(1+k),\ldots,k/(1+k), \underbrace{1,\ldots,1}_{q \text{ times}}).$$

*Define*

$$\Theta_j = \{\theta \in R^{m_j} : \|D^r f_{\theta,j}\|_{L_\infty[0,1]} \leq L, \ \forall \, r \in \{0,1,\ldots,q-1\} \text{ and } \|f_{\theta,j}\|_\infty \leq M\},$$

where for $\theta = (\theta_1, \ldots, \theta_{m_j}) \in R^{m_j}$,

$$f_{\theta,j} = \sum_{i=1}^{m_j} \theta_i B_{j,i} \stackrel{def}{=} \theta' B. \tag{21}$$

Define $\eta_j$ according to (18) with

$$A_j = 9.64\sqrt{q}(2q+1)9^{q-1} + 0.06 \quad \text{and} \quad C_j = m_j + L, \tag{22}$$

and define $a_j$ according to (17). Let $\pi_j$ to be the Lebesgue measure on $\Theta_j$. Let $\tilde{\pi}_j$ be the induced prior of $\pi_j$ and $\tilde{B}_{L_2(\mu)}(s_n)$ denote the $s_n$ $L_2(\mu)$ neighborhood of $f_o$, as defined on page 1557. Then for the prior $\tilde{\pi} = \sum_j a_j \tilde{\pi}_j$, the posterior probability $\tilde{\pi}(\tilde{B}_{L_2(\mu)}(s_n)^c | X_1, \ldots, X_n)$ converges to zero in probability for some $s_n \propto n^{-s/(1+2s)}$.

The proof for Lemma 5 is given in Section 4.

Here is a lemma that is useful for verifying Assumption 7 to prove Lemma 5.

LEMMA 6.  *Suppose that $\{S_j : j \in J\}$ is a countable collection of linear function spaces on $[0,1]$. Suppose that for each $S_j$ there is a basis $\{B_{j,1}, \ldots, B_{j,m_j}\}$. Suppose that there exist constants $T_1$ and $T_2$ such that for $\theta = (\theta_1, \ldots, \theta_{m_j}) \in R^{m_j}$,*

$$\left\| \sum_{i=1}^{m_j} \theta_i B_{j,i} \right\|_\infty \leq T_1 \max_i |\theta_i| \tag{23}$$

*and*

$$\left\| \sum_{i=1}^{m_j} \theta_i B_{j,i} \right\|_2 \geq \frac{T_2}{\sqrt{m_j}} \sqrt{\sum_{i=1}^{m_j} \theta_i^2}, \tag{24}$$

*where $\|\cdot\|_2$ denotes the $L_2$ norm with respect to the Lebesgue measure on $[0,1]$. Suppose that for $j \in J$, $\Theta_j \subset R^{m_j}$ and $f_{\theta,j}$ is as defined in (21). Then Assumption 7 holds with*

$$A_j = 9.64 \frac{T_1}{T_2} + 0.06. \tag{25}$$

The proof is a straightforward modification of the proof for Lemma 1 of Yang and Barron (1998).

**4. Proofs.**

4.1. *Proof of Theorem* 1. We prove Theorem 1 by giving bounds for $U_n$ and $V_n$, respectively, and then combining the bounds to show that $U_n/V_n$ converges to zero. For finding an upper bound for $U_n$, we would like to use the following lemma, which is a modified version of Lemma 0 by Yang and Barron (1998).

LEMMA 7. *Suppose that Assumption* 1 *holds and*
$$\frac{\xi_j}{m_j} \geq \frac{4}{1-4\gamma} \log\left(\frac{46.2 A_j \sqrt{1-4\gamma}}{\gamma}\right).$$

*Then*
$$P_o^*\left[ \text{for some } \theta \in \Theta_j, \frac{1}{n}\sum_{i=1}^n \log \frac{f_{\theta,j}(X_i)}{f_o(X_i)} \geq -\gamma d_H^2(f_o, f_{\theta,j}) + \frac{\xi_j}{n} \right]$$
$$\leq 15.1 \exp\left(-\frac{1-4\gamma}{8}\xi_j\right),$$

*where $P_o^*$ is the outer measure for $P_{f_o}^n$.*

PROOF. Suppose that Assumption 1 holds. We will show that for any $r > 0$ and $\delta \leq 0.056r$,

(26) $$N(B_{d_H,j}(r), \delta, d_{j,\infty}) \leq \left(\frac{3 A_j r}{\delta}\right)^{m_j},$$

where $B_{d_H,j}(r)$ is as defined on page 1557. Then the result in Lemma 7 follows from Lemma 0 in Yang and Barron (1998).

Below is the proof of (26). Fix $\varepsilon > 0$. Let $\theta_* \in \Theta_j$ be such that
$$d_H(f_o, f_{\theta_*,j}) \leq \varepsilon r + \inf_{\theta \in \Theta_j} d_H(f_o, f_{\theta,j}).$$

Then for $\theta \in \Theta_j$,
$$d_H(f_o, f_{\theta,j}) \geq \frac{1}{2}(d_H(f_o, f_{\theta_*,j}) + d_H(f_o, f_{\theta,j})) - \frac{\varepsilon r}{2}$$
$$\geq \frac{1}{2} d_H(f_{\theta,j}, f_{\theta_*,j}) - \frac{\varepsilon r}{2},$$

so we have
$$B_{d_H,j}(r) = \{\theta \in \Theta_j : d_H(f_o, f_{\theta,j}) \leq r\}$$
$$\subset \{\theta \in \Theta_j : d_H(f_{\theta,j}, f_{\theta_*,j}) \leq (2+\varepsilon)r\}$$
$$= B_{d_H,j}(\theta_*, (2+\varepsilon)r),$$

where $B_{d_H,j}(\theta_*, (2+\varepsilon)r)$ is as defined on page 1558. Take $\varepsilon = 1$; then by Assumption 1, for any $r > 0$ and $\delta \leq 0.056r$, (26) holds, so by Lemma 0 in Yang and Barron (1998) the proof for Lemma 7 is complete. □

Suppose Assumption 1 holds. Let $a_j$ and $\eta_j$ be as specified in (2) and (3) take $\xi_j = \eta_j + \gamma n s_n^2 / 2$. Then by Lemma 7 and we have

$$U_n \leq \left(\sum_j a_j e^{\xi_j}\right) e^{-\gamma n s_n^2}$$

$$= \alpha e^{-\gamma n s_n^2/2} \sum_j \exp\left(-\frac{1-4\gamma}{8}\eta_j\right) \leq \alpha e^{-\gamma n s_n^2/2}$$

except on a set of probability no greater than

$$\sum_j 15.1 \exp\left(-\frac{1-4\gamma}{8}\xi_j\right)$$

$$= 15.1 \exp\left(-\frac{(1-4\gamma)\gamma n s_n^2}{16}\right) \sum_j \exp\left(-\frac{1-4\gamma}{8}\eta_j\right)$$

$$\leq 15.1 \exp\left(-\frac{(1-4\gamma)\gamma n s_n^2}{16}\right).$$

That is, an upper bound for $U_n$ is given by

(27) $$P_{f_o}^n[U_n > \alpha e^{-\gamma n s_n^2/2}] \leq 15.1 \exp\left(-\frac{(1-4\gamma)\gamma n s_n^2}{16}\right).$$

To find a lower bound for $V_n$, we will use Lemma 1 of Shen and Wasserman (2001). Let

$$\tilde{B}_D(r) = \{g : D(f_o \| g) \leq r, V'(f_o \| g) \leq r\},$$

where $V'(f \| g) = \int f(\log(f/g) - D(f \| g))^2 \, d\mu$. Here is the lemma.

LEMMA 8. *For $t_n > 0$,*

$$P_{f_o}^n\left(V_n \leq \frac{1}{2}\tilde{\pi}(\tilde{B}_D(t_n))e^{-2nt_n}\right) \leq \frac{2}{nt_n}.$$

Suppose that Assumptions 2–5 hold. Let $B_{d_{j_n},j_n}(\theta, \varepsilon_n)$ denote the $d_{j_n}$-ball centered at $\theta$ with radius $\varepsilon_n$ in $\Theta_{j_n}$ and define

$$B_{D,j_n}(t_n) = \{\theta \in \Theta_{j_n} : D(f_o \| f_{\theta,j_n}) \leq t_n, V(f_o \| f_{\theta,j_n}) \leq t_n\}.$$

We will first show that

(28) $$B_{d_{j_n},j_n}(\beta_n, \varepsilon_n) \subset B_{D,j_n}(t_n)$$

for some $t_n \propto \varepsilon_n^2$ and that

(29) $$\pi_{j_n}(B_{d_{j_n},j_n}(\beta_n, \varepsilon_n)) \geq \left(\frac{1}{A_{j_n}^{b_1+b_2} K_4 K_5}\right)^{m_{j_n}}.$$

Then we will deduce a lower bound for $\tilde{\pi}(\tilde{B}_D(t_n))$ based on (28) and (29) to apply Lemma 8.

To prove (28), note that for $\theta \in B_{d_{j_n},j_n}(\beta_n,\varepsilon_n)$, by Assumptions 2 and 3 we have

$$V(f_o\|f_{\theta,j_n}) \leq 2\varepsilon_n^2 + 2K_0'\varepsilon_n^2$$

and

$$D(f_o\|f_{\theta,j_n}) \leq K_0'' V(f_o\|f_{\theta,j_n}) \leq 2K_0''(1+K_0')\varepsilon_n^2.$$

Therefore, (28) holds for $t_n = 2\max(1,K_0'')(1+K_0')\varepsilon_n^2 \stackrel{\text{def}}{=} K'\varepsilon_n^2$.

To prove (29), note that by Assumption 4 there exist $\theta_1,\ldots,\theta_{d^*} \in \Theta_{j_n}$ such that

$$d^* \leq (A_{j_n}^{b_1} K_4)^{m_{j_n}} \quad \text{and} \quad \bigcup_{i=1}^{d^*} B_{d_{j_n},j_n}(\theta_i,\varepsilon_n) \supset \Theta_{j_n},$$

so

$$\pi_{j_n}(B_{d_{j_n},j_n}(\beta_n,\varepsilon_n)) \geq \frac{\pi_{j_n}(B_{d_{j_n},j_n}(\beta_n,\varepsilon_n))}{\sum_{i=1}^{d^*} \pi_{j_n}(B_{d_{j_n},j_n}(\theta_i,\varepsilon_n))}$$
$$\geq \left(\frac{1}{A_{j_n}^{b_1+b_2} K_4 K_5}\right)^{m_{j_n}},$$

where the last inequality follows from Assumption 5.

It is clear that

$$\tilde{\pi}(\tilde{B}_D(t_n)) \geq a_{j_n}\pi_{j_n}(B_{D,j_n}(t_n))$$
$$\stackrel{(28)}{\geq} a_{j_n} B_{d_{j_n},j_n}(\beta_n,\varepsilon_n)$$
$$\stackrel{(29)}{\geq} a_{j_n}\left(\frac{1}{A_{j_n}^{b_1+b_2} K_4 K_5}\right)^{m_{j_n}},$$

so by Lemma 8, we have that except on a set of probability no greater than $2/(nt_n)$,

$$V_n \geq \frac{1}{2}e^{-2nt_n} a_{j_n}\pi_{j_n}(B_{D,j_n}(t_n))$$
$$\geq \frac{e^{-2nt_n}}{2}\alpha\exp\left(-\left(1+\frac{1-4\gamma}{8}\right)\eta_{j_n}\right)\left(\frac{1}{A_{j_n}^{b_1+b_2} K_4 K_5}\right)^{m_{j_n}}$$
(30) $$\geq \frac{\alpha}{2}\exp\left(-2nt_n - \eta_{j_n}\left(1+\frac{1-4\gamma}{8}+b_1+b_2+(\log(K_4 K_5))_+\right)\right)$$

$$\overset{(4)}{\geq} \frac{\alpha}{2} \exp\left(-2nt_n - n\varepsilon_n^2\left(1 + \frac{1-4\gamma}{8} + b_1 + b_2 + (\log(K_4 K_5))_+\right)\right)$$

$$= \frac{\alpha}{2} e^{-Kn\varepsilon_n^2},$$

where $K = 2K' + 1 + (1-4\gamma)/8 + b_1 + b_2 + (\log(K_4 K_5))_+$. Here the third inequality follows from the fact that

$$\frac{\eta_j}{m_j} \geq \frac{4}{1-4\gamma} \log\left(\frac{46.2 A_j \sqrt{1-4\gamma}}{\gamma}\right) \overset{A_j \geq 0.0056 = 0.13\gamma/\sqrt{1-4\gamma}}{\geq} \max(1, \log A_j)$$

for all $j$.

Now we will bound $U_n/V_n$ by combining (27) and (30). In (27) set $s_n^2 = 4K\varepsilon_n^2/\gamma$. Then

$$\tilde{\pi}(\tilde{B}_{d_\mathrm{H}}(s_n)^c | X_1, \ldots, X_n) = \frac{U_n}{V_n} \leq 2\exp(-Kn\varepsilon_n^2)$$

except on a set of probability no greater than

$$15.1 \exp\left(-\frac{(1-4\gamma)Kn\varepsilon_n^2}{4}\right) + \frac{2}{K' n \varepsilon_n^2},$$

which converges to zero because $n\varepsilon_n^2 \to \infty$ by Assumption 6.

4.2. *Proof of Lemma* 1. We will verify Assumptions 1–6 for the spline example. To verify Assumption 1, we will apply Lemma 2. From page 143 (4.80) in Schumaker (1981)

$$\left\|\sum_{i=1}^{m_j} \theta_i B_{j,i}\right\|_\infty \leq \max_i |\theta_i|.$$

Since $m_j$ and $B_{j,i}$ depend on $(k,q)$ but not on $L$, we set $l = (k,q)$, $m_l = m_j$ and $B_{l,i} = B_{j,i}$. Then (8) holds with $T_1 = 1$. To check (9), note that from (4.79) and (4.86) in Schumaker (1981), we have that for each $i \in \{1, \ldots, m_l\}$,

$$|\theta_i| \leq (2q+1)9^{q-1}(y_{i+q} - y_i)^{-1/2} \|\theta_i B_{l,i}\|_{L_2[y_i, y_{i+q}]},$$

where $y_1, \ldots, y_{2q+k}$ are as defined in Lemma 1 and $L_2[y_i, y_{i+q}]$ is the $L_2$ metric with respect to the Lebesgue measure on $[y_i, y_{i+q}]$. Since $y_{i+q} - y_i \geq 1/(1+k)$,

$$\sum_{i=1}^{m_l} \theta_i^2 \leq (2q+1)^2 9^{2(q-1)} (k+1) \sum_{i=1}^{m_l} \|\theta_i B_{l,i}\|_{L_2[y_i, y_{i+q}]}^2$$

$$\leq (2q+1)^2 9^{2(q-1)} (k+q) q \left\|\sum_{i=1}^{m_l} \theta_i B_{l,i}\right\|_2^2,$$

which implies that (9) holds with $T_2 = 1/(\sqrt{q}(2q+1)9^{q-1})$. By Lemma 2, Assumption 1 holds for $A_j$ and $m_j$ in (7). Also note that for the $C_j$ specified in (7), $\sum_j e^{-C_j} = e^{-2}/(1-e^{-1})^3 < 1$ as required.

To verify Assumption 2, we need to find $j_n$ and $\beta_n$. Take $j_n = (k_n, q^*, L^*)$, where $\{k_n\}$ is a sequence of positive integers such that

$$c_3 n^{1/(1+2s)} \le k_n \le c_4 n^{1/(1+2s)} \qquad \text{for all } n$$

for some constants $c_3$ and $c_4$, $q^* = s + 1$, and

$$L^* = \min\{L : L \text{ is a positive integer}, L \ge 2^s + \alpha_{q^*} M_0 + M_0\},$$

where $M_0 = \max_{0 \le r \le s} \|D^r \log f_o\|_{L_\infty}$. To control the error $\max(D(f_o\|f_{\beta_n, j_n}), V(f_o\|f_{\beta_n, j_n}))$, we use the following fact.

FACT 1. For $j$ such that $q \ge s + 1$, there exists $\beta \in R^{m_j}$ such that

(31)
$$\|D^r(\log f_o - \log f_{\beta,j})\|_\infty \le \alpha_q \left(\frac{1}{k+1}\right)^{s-r} M_0 \qquad \text{for } 0 \le r \le s-1,$$
$$\|D^s \log f_{\beta,j}\|_\infty \le \alpha_q M_0.$$

This fact follows from (6.50) in Schumaker (1981) and the result that for $\theta = (\theta_1, \ldots, \theta_{m_j}) \in R^{m_j}$,

$$|\psi(\theta)| = \left|\log \int_0^1 \exp\left(-\log f_o(x) + \sum_{i=1}^{m_j} \theta_i B_{j,i}(x)\right) f_o(x)\, dx\right|$$
$$\le \left\|\log f_o - \sum_{i=1}^{m_j} \theta_i B_{j,i}\right\|_\infty.$$

From the fact, there exists $\beta_n \in R^{m_{j_n}}$ such that

$$\|\log f_o - \log f_{\beta_n, j_n}\|_\infty \le \alpha_{q^*} M_0 \left(\frac{1}{k_n + 1}\right)^s.$$

Since $D(f_o\|f_{\beta_n, j_n})$ and $V(f_o\|f_{\beta_n, j_n})$ are bounded by $\|\log f_o - \log f_{\beta_n, j_n}\|_\infty$, we have

$$\max(D(f_o\|f_{\beta_n, j_n}), V(f_o\|f_{\beta_n, j_n})) + \frac{\eta_{j_n}}{n}$$
$$\le \alpha_{q^*} M_0 \left(\frac{1}{k_n + 1}\right)^{2s} + \frac{c_2 k_n}{n} \le c_1 n^{-2s/(1+2s)}$$

for some constants $c_1$ and $c_2$. So Assumption (2) holds if $\beta_n \in \Theta_{j_n}$ and

(32) $$\varepsilon_n^2 = c_1 n^{-2s/(1+2s)}.$$

To verify that $\beta_n \in \Theta_{j_n}$, we need to make sure $\beta_n' \mathbb{1}_{m_{j_n}} = 0$ and $\max_{0 \leq r \leq q-1} \|D^r \log f_{\beta_n, j_n}\|_{L_\infty} \leq L^*$. For the first condition, $\beta_n' \mathbb{1}_{m_{j_n}} = 0$, we can assume it without loss of generality, because $\log f_{\beta_n, j_n}$ does not change when $\beta_n$ is shifted by a constant. The second condition holds because of the second equation in (31).

Now let us verify Assumptions 3–5 with $d_{j_n} = d_{j_n, \infty}$, where $d_{j_n, \infty}$ is as defined in Assumption 1. For the verification of Assumption 3, we will use the following fact.

FACT 2. Suppose that

$$(33) \qquad \int f_o \left( \log \frac{f_{\eta, j_n}}{f_{\theta, j_n}} \right)^2 \leq K_0 d_{j_n}^2(\eta, \theta) \qquad \text{for all } \eta, \theta \in \Theta_{j_n}$$

for some constant $K_0$ and

$$(34) \qquad \sup_{\theta \in \Theta_{j_n}} \|\log f_o - \log f_{\theta, j_n}\|_\infty \leq \log K_3$$

for some constant $K_3$. Then Assumption 3 holds with $K_0' = K_0$ and $K_0'' = K_3/2$.

The proof of the fact is a straightforward application of an equation in Lemma 1 by Barron and Sheu (1991), which gives

$$(35) \qquad D(f_o \| f_{\theta, j_n}) \leq \tfrac{1}{2} e^{\|\log f_o - \log f_{\theta, j_n}\|_\infty} V(f_o \| f_{\theta, j_n})$$

for all $\theta \in R^{m_{j_n}}$. It is clear that (33) holds with $K_0 = 1$ and that (34) holds with $K_3 = e^{2L^*}$, so by Fact 2, Assumption 3 holds.

For Assumption 4, by Theorems IV and XIV of Kolmogorov and Tikhomirov (1961), there exists an $\varepsilon_n$-net $F_{\varepsilon_n}$ for $\Theta_{j_n}$ with respect to $d_{j_n}$ so that

$$\log \operatorname{card}(F_{\varepsilon_n}) \leq c_{q^*, L^*} \left( \frac{1}{\varepsilon_n} \right)^{1/(q^*-1)}$$

$$= c_{q^*, L^*} \left( \frac{1}{\varepsilon_n} \right)^{1/s} \leq c_{q^*, L^*}(k_n + 1) \leq c_{q^*, L^*} m_{j_n}.$$

Therefore, Assumption 4 holds with $K_4 = e^{c_{q^*, L^*}}$ and $b_1 = 0$.

We will check Assumption 5. For a positive integer $m$, for $t = (t_1, \ldots, t_m) \in R^m$, define

$$\|t\|_\infty = \max_{1 \leq i \leq m} |t_i|.$$

To bound $\pi_{j_n}(B_{d_{j_n}, j_n}(\beta_n, \varepsilon_n))$, we will show that

$$(36) \quad \left\{ \theta \in R^{m_{j_n}} : \theta' \mathbb{1}_{m_{j_n}} = 0, \|\theta - \beta_n\|_\infty \leq c_6 \left( \frac{1}{k_n + 1} \right)^s \right\} \subset B_{d_{j_n}, j_n}(\beta_n, \varepsilon_n),$$

where $c_6 = \min(1, \sqrt{c_1}/2(\sup_n n^{s/(1+2s)}(k_n+1)^{-s}))$. To prove (36), suppose that $\theta \in R^{m_{j_n}}$ and

$$\theta' \mathbb{1}_{m_{j_n}} = 0 \quad \text{and} \quad \|\theta - \beta_n\|_\infty \leq c_6 \left(\frac{1}{k_n+1}\right)^s.$$

We will show that

(37) $$d_{j_n}(\theta, \beta_n) \leq \varepsilon_n$$

and

(38) $$\theta \in \Theta_{j_n}.$$

Inequality (37) holds since

$$\|\log f_{\theta,j_n} - \log f_{\beta_n,j_n}\|_\infty \leq 2\|\theta - \beta_n\|_\infty \leq 2c_6 c_5 n^{-s/(1+2s)} \leq \varepsilon_n,$$

where $c_5 = \sup_n (k_n+1)^{-s} n^{s/(1+2s)}$. Here the second inequality holds because

$$|\psi(\theta) - \psi(\beta_n)| = \left|\log \int e^{(\theta-\beta_n)'B} e^{\beta_n' B - \psi(\beta_n)}\right| \leq \|\theta - \beta_n\|_\infty.$$

To prove (38), we need the following inequality:

(39) $$\|D^r(\theta'B - \beta'B)\|_{L_\infty} \leq 2^r (k+1)^r \|\theta - \beta\|_\infty \qquad \text{for all } 0 \leq r \leq s,$$

which is deduced from (4.54) in Schumaker (1981). Now note that for $0 < r < s$,

$$\begin{aligned}
\|D^r \log f_{\theta,j_n}\|_\infty &= \|D^r \theta' B\|_\infty \\
&\leq \|D^r(\theta'B - \beta_n' B)\|_\infty + \|D^r(\beta_n' B - \log f_o)\|_\infty \\
&\quad + \|D^r \log f_o\|_\infty \\
&\stackrel{(39),(31)}{\leq} 2^r (k_n+1)^r \|\theta - \beta_n\|_\infty + \alpha_{q^*} M_0 \left(\frac{1}{k_n+1}\right)^{s-r} + M_0 \\
&\leq \left(\frac{1}{k_n+1}\right)^{s-r} (2^r + \alpha_{q^*} M_0) + M_0 \leq L^*,
\end{aligned}$$

for $r = 0$,

$$\begin{aligned}
\|\log f_{\theta,j_n}\|_\infty &\leq \|\log f_{\theta,j_n} - \log f_{\beta_n,j_n}\|_\infty + \|\log f_{\beta_n,j_n} - \log f_o\|_\infty + \|\log f_o\|_\infty \\
&\leq 2\|\theta - \beta_n\|_\infty + \|\log f_{\beta_n,j_n} - \log f_o\|_\infty + M_0 \\
&\leq \left(\frac{1}{k_n+1}\right)^s (2 + \alpha_{q^*} M_0) + M_0 \leq L^*,
\end{aligned}$$

and for $r = s$,

$$\|D^s \log f_{\theta,j_n}\|_{L_\infty} = \|D^s \theta' B\|_{L_\infty}$$
$$\leq \|D^s(\theta' B - \beta'_n B)\|_{L_\infty} + \|D^s \beta'_n B\|_{L_\infty}$$
$$\overset{(39),(31)}{\leq} 2^s + \alpha_{q^*} M_0 \leq L^*.$$

Therefore, $\theta \in \Theta_{k_n,q^*,L^*}$, so (38) and (36) hold. To bound $\pi_{j_n}(B_{d_{j_n},j_n}(\theta_1,\varepsilon_n))$ in Assumption 5, note that for all $\varepsilon > 0$ and for all $j$,

(40) $\quad \{\theta \in \Theta_j : \|\log f_{\theta,j} - \log f_{\theta_1,j}\|_\infty \leq \varepsilon\} \subset \{\theta \in \Theta_j : \|\theta - \theta_1\|_\infty \leq 2\beta^*_{q^*}\varepsilon\},$

where $\beta^*_{q^*}$ is some positive constant. This result follows from Lemma 4.3 of Ghosal, Ghosh and van der Vaart (2000), which implies that for all $\theta$, $\theta_1 \in R^{m_{j_n}}$,

$\|\theta - \theta_1\|_\infty \leq \|\log f_{(\theta-\theta_1),j_n}\|_\infty$ times some constant depending on $q^*$,

and from the fact that

$$\|\log f_{(\theta-\theta_1),j_n} - (\log f_{\theta_1,j_n} - \log f_{\theta,j_n})\|_\infty$$
$$= |\psi(\theta - \theta_1) - (\psi(\theta) - \psi(\theta_1))|$$
$$= \left|\log \int \exp(\theta' B - \psi(\theta) - (\theta'_1 B - \psi(\theta_1)))\right|$$
$$\leq \|\log f_{\theta_1,j_n} - \log f_{\theta,j_n}\|_\infty.$$

Then by (40) and by (36) we have

$$\frac{\pi_{j_n}(B_{d_{j_n},j_n}(\beta_n,\varepsilon_n))}{\pi_{j_n}(B_{d_{j_n},j_n}(\theta_1,\varepsilon_n))} \geq \frac{(c_6(1/(k_n+1))^s)^{k_n+q^*-1}}{(\beta^*_{q^*}\varepsilon_n)^{k_n+q^*-1}}$$
$$\geq \left(\frac{c_6}{\beta^*_{q^*}\varepsilon_n(1+(c_4\sqrt{c_1}/\varepsilon_n)^{1/s})^s}\right)^{k_n+q^*-1}.$$

For $n$ such that $0 < \varepsilon_n \leq 1$,

$$\frac{\pi_{j_n}(B_{d_{j_n},j_n}(\beta_n,\varepsilon_n))}{\pi_{j_n}(B_{d_{j_n},j_n}(\theta_1,\varepsilon_n))} \geq \left(\frac{c_6}{\beta^*_{q^*}\varepsilon_n((1/\varepsilon_n)^{1/s}+(c_4\sqrt{c_1}/\varepsilon_n)^{1/s})^s}\right)^{k_n+q^*-1}$$
$$= \left(\frac{c_6}{\beta^*_{q^*}(1+(c_4\sqrt{c_1})^{1/s})^s}\right)^{k_n+q^*-1}.$$

Without loss of generality, we can assume that $\beta^*_{q^*} > 1$, so it is clear that Assumption 5 holds with $K_5 = \beta^*_{q^*}(1+(c_4\sqrt{c_1})^{1/s})^s/c_6$ and $b_2 = 0$.

For Assumption 6, it should be clear that it holds with the $\varepsilon_n$ specified in (32). Now by Theorem 1, the result in Lemma 1 holds.

4.3. *Proof of Lemma* 3. We will verify Assumptions 1–6 for the Haar basis example. To verify Assumption 1, we will apply Lemma 4. First, by (3.7) in Barron, Birgé and Massart (1999), (14) holds for $K_l = 2^{(l+1)/2}$. Second, for all $j$ and $\theta \in \Theta_j$, $|\phi(\theta)| = |\log \int e^{\theta' B}| \leq \|\theta' B\|_\infty$, so (15) holds. Therefore, by Lemma 4, Assumption 1 holds for $A_j$ and $m_j$ in (13). Note that for the $C_j$ specified in (13), $\sum_j e^{-C_j} < 1$ as required.

To verify Assumption 2, we will first choose $j_n$ and $\beta_n$, and then show that

$$\|\log f_o - \log f_{\beta_n, j_n}\|_2 \leq c_{1,\alpha,f_o,H_0} \left(\frac{1}{m_{j_n}}\right)^\alpha,$$

(41)

$$\|\log f_o - \log f_{\beta_n, j_n}\|_\infty \leq 2c_{2,f_o}$$

for some constants $c_{1,\alpha,f_o,H_0}$ and $c_{2,f_o}$ and that $\beta_n \in \Theta_{j_n}$. Then we will take $\varepsilon_n$ according to an upper bound for the left-hand side of (31) so that Assumption 2 holds. We will see that $\varepsilon_n$ converges to zero at the rate $(\log n)^{1/2} n^{-\alpha/(1+2\alpha)}$ as required.

$j_n$ and $\beta_n$ are defined as follows. Let $\{l_n\}$ be a sequence of integers such that

$$k_3 n^{1/(1+2\alpha)} \leq 2^{l_n+1} \leq k_4 n^{1/(1+2\alpha)},$$

where $k_3$ and $k_4$ are positive constants. Let

$$\beta_0 + \sum_{i=1}^{m_{j_n}-1} \beta_{l_n,i} B_{l_n,i} \stackrel{\text{def}}{=} \beta_0 + \beta_n' B$$

be the $L_2$ projection of $\log f_o$ to the space spanned by 1 and $B_{l_n,i} : i = 1, \ldots, m_{j_n} - 1$. Let $M_0 = \|\log f_o\|_\infty$ and $c_{2,f_o} = \sup_n \|\log f_o - \beta_0 - \beta_n' B\|_\infty$. ($c_{2,f_o}$ is finite since $\beta_0 + \beta_n' B$ converges to $\log f_o$ uniformly.) Define

$$L^* = \min\{L : L \text{ is a positive integer and } L \geq 2c_{2,f_o} + 3M_0\}.$$

Set $j_n = (l_n, L^*)$.

To prove (41), we will bound $\log f_o - \beta_0 - \beta_n' B$ and $\beta_0 + \psi(\beta_n)$, respectively. By (12) we have

$$\|\log f_o - \beta_0 - \beta_n' B\|_2 \leq \frac{H_0 2^{-\alpha(l_n+1)}}{\sqrt{1-2^{-2\alpha}}} \leq \frac{H_0}{\sqrt{1-2^{-2\alpha}}} \left(\frac{1}{m_{j_n}}\right)^\alpha.$$

To bound $\beta_0 + \psi(\beta_n)$, let $\Delta = \int (e^{\beta_0 + \beta_n' B - \log f_o} - 1) f_o$ and $b = \|\log f_o - \beta_0 - \beta_n' B\|_\infty$. Then

$$|\beta_0 + \psi(\beta_n)| = \left|\log \int e^{\beta_0 + \beta_n' B - \log f_o} f_o\right|$$

$$= |\log(1 + \Delta)|$$

$$\leq \max\left(\Delta, \frac{-\Delta}{1+\Delta}\right)$$

$$\leq |\Delta| e^{b+M_0} (\text{since } e^{-b-M_0} \leq 1 + \Delta \leq e^{b+M_0})$$

$$\leq e^{b+2M_0}\left(1 + \frac{1}{2}e^b \|\log f_o - \beta_0 - \beta'_n B\|_2\right) \|\log f_o - \beta_0 - \beta'_n B\|_2,$$

where the last inequality follows from the Cauchy–Schwarz inequality and (3.3) in Barron and Sheu (1991), which says that

$$\frac{z^2}{2} e^{-\max(-z,0)} \leq e^z - 1 - z \leq \frac{z^2}{2} e^{\max(z,0)} \qquad \text{for all } z.$$

Therefore, the first inequality in (41) holds. The second inequality in (41) also holds since

$$\|\log f_o - \beta_0 - \beta'_n B\|_\infty \leq \|\log f_o - \beta_0 - \beta'_n B\|_\infty + |\beta_0 + \psi(\beta_n)|$$

$$= c_{2,f_o} + \left|\log \int e^{\beta_0 + \beta'_n B - \log f_o} f_o\right| \leq 2c_{2,f_o}.$$

Now we have proved (41), which implies that $\|\log f_{\beta_n, j_n}\|_\infty \leq L^*$, so $\beta_n \in \Theta_{j_n}$.

The $L_2$ bound in (41) gives a bound for the error $\max(D(f_o\|f_{\beta_n,j_n}), V(f_o\|f_{\beta_n,j_n}))$ since

$$(42) \quad V(f_o\|f_{\beta_n,j_n}) = \int f_o \left(\log \frac{f_o}{f_{\beta_n,j_n}}\right)^2 \leq e^{\|\log f_o\|_\infty} \|\log f_o - \log f_{\beta_n,j_n}\|_2^2$$

and by (35) and (41),

$$(43) \qquad D(f_o\|f_{\beta_n,j_n}) \leq \tfrac{1}{2} e^{2c_{2,f_o}} V(f_o\|f_{\beta_n,j_n}).$$

By (41)–(43) and the definition of $\eta_{j_n}$, we can find two constants $k_1$ and $k_2$ which depend only on $\alpha$, $f_o$ and $H_0$ such that

$$\max(D(f_o\|f_{\beta_n,j_n}), V(f_o\|f_{\beta_n,j_n})) + \frac{\eta_{j_n}}{n} \leq k_1\left(\frac{1}{m_{j_n}}\right)^{2\alpha} + k_2\frac{m_{j_n} \log m_{j_n}}{n}.$$

Since $l_n$ is chosen such that $k_3 n^{1/(1+2\alpha)} \leq m_{j_n} \leq k_4 n^{1/(1+2\alpha)}$, we have

$$\max(D(f_o\|f_{\beta_n,j_n}), V(f_o\|f_{\beta_n,j_n})) + \frac{\eta_{j_n}}{n}$$

$$\leq \left(\frac{k_1}{k_3^{2\alpha}} + k_2 k_4 \log k_4 + \frac{k_2 k_4}{1+2\alpha}\right) n^{-2\alpha/(1+2\alpha)} \log n$$

$$\stackrel{\text{def}}{=} k_5 n^{-2\alpha/(1+2\alpha)} \log n.$$

Hence, Assumption 2 holds with $\varepsilon_n^2 = k_5 n^{-2\alpha/(1+2\alpha)} \log n$.

To verify Assumption 3, for all positive integers $m$ and for all $t = (t_1, \ldots, t_m) \in R^m$ define

$$\|t\| = \sqrt{\sum_{i=1}^m t_i^2}.$$

Let $d_{j_n} = \|\cdot\|$ on $R^{m_{j_n}-1}$. We will verify Assumption 3 using Fact 2. For $\eta, \theta \in \Theta_{j_n}$, since

$$\psi(\eta) - \psi(\theta) = \log \int e^{(\theta-\eta)'B} f_{\eta,j_n}$$

$$\leq \log \int (1 + (\theta-\eta)'B e^{(\theta-\eta)'B}) f_{\eta,j_n}$$

$$\leq \log\left(1 + \sqrt{\int ((\theta-\eta)'B)^2} \sqrt{\int e^{2(\theta-\eta)'B} f_{\eta,j_n}^2}\right)$$

$$\leq \log(1 + \|\theta - \eta\| e^{4L^*})$$

$$\leq e^{4L^*} \|\theta - \eta\|,$$

$$\|\log f_{\eta,j_n} - \log f_{\theta,j_n}\|_2^2 = (\psi(\eta) - \psi(\theta))^2 + \|\eta - \theta\|^2$$

and

$$\int f_o \left(\log \frac{f_{\eta,j_n}}{f_{\theta,j_n}}\right)^2 \leq e^{\|\log f_o\|_\infty} \|\log f_{\eta,j_n} - \log f_{\theta,j_n}\|_2^2$$

$$= e^{M_0} \|\log f_{\eta,j_n} - \log f_{\theta,j_n}\|_2^2,$$

(33) holds with $K_0 = e^{M_0}(1 + e^{8L^*})$ and clearly, (34) holds with $K_3 = e^{M_0 + 2L^*}$. Therefore, by Fact 2, Assumption 3 holds.

For checking Assumption 4, note that

$$\Theta_{j_n} \subset \{\theta \in R^{m_{j_n}-1} : \|\theta\|_\infty \leq L^*\},$$

which implies that for every $\varepsilon > 0$, there exists an $\varepsilon$-net $F_\varepsilon$ for $\Theta_{j_n}$ with respect to $\|\cdot\|_\infty$ so that

$$\mathrm{card}(F_{\varepsilon_n}) \leq \left(1 + \frac{2L^*}{\varepsilon}\right)^{m_{j_n}-1}.$$

By the fact that $\|\theta\| \leq \sqrt{m_{j_n} - 1} \|\theta\|_\infty$ for all $\theta \in \Theta_{j_n}$, there exists an $\varepsilon_n$-net $F_{\varepsilon_n}$ for $\Theta_{j_n}$ with respect to $d_{j_n}$ such that

$$\mathrm{card}(F_{\varepsilon_n}) \leq \left(1 + \frac{2L^* \sqrt{m_{j_n} - 1}}{\varepsilon_n}\right)^{m_{j_n}-1}.$$

Since
$$\frac{1 + (2L^*\sqrt{m_{j_n} - 1})/\varepsilon_n}{A_{j_n}^{3\alpha}} \leq \frac{(1 + 2L^*\sqrt{k_4}/k_5)n^{1.5\alpha/(1+2\alpha)}}{k_3^{1.5\alpha}n^{1.5\alpha/(1+2\alpha)}},$$

Assumption 4 holds with $K_4 = (1 + 2L^*\sqrt{k_4}/k_5)/(k_3^{1.5\alpha})$ and $b_1 = 3\alpha$.

For Assumption 5, to bound $\pi_{j_n}(B_{d_{j_n},j_n}(\beta_n, \varepsilon_n))$, we will show that

(44) $\quad \left\{\theta \in R^{m_{j_n}-1} : \|\theta - \beta_n\|_\infty \leq \dfrac{\varepsilon_n}{m_{j_n}\sqrt{m_{j_n} - 1}}\right\} \subset B_{d_{j_n},j_n}(\beta_n, \varepsilon_n)$

for $n$ such that $\varepsilon_n \leq M_0$. For $\theta \in R^{m_{j_n}-1}$ such that $\|\theta - \beta_n\|_\infty \leq \varepsilon_n/(m_{j_n}\sqrt{m_{j_n} - 1})$,

$$\|\theta - \beta_n\| \leq \sqrt{m_{j_n} - 1}\|\theta - \beta_n\|_\infty \leq \frac{\varepsilon_n}{m_{j_n}} \leq \varepsilon_n,$$

so it suffices to show that $\theta \in \Theta_{j_n}$. For $n$ such that $\varepsilon_n \leq M_0$,

$$\|\theta'B\|_\infty \leq \|\theta'B - \beta_n'B\|_\infty + \|\beta_0 + \beta_n'B - \log f_o\|_\infty + |\beta_0| + \|\log f_o\|_\infty$$
$$\leq m_{j_n}\|\theta - \beta_n\| + 2c_{2,f_o}M_0 + 2M_0$$
$$\leq \varepsilon_n + 2c_{2,f_o}M_0 + 2M_0$$
$$\leq 2c_{2,f_o}M_0 + 3M_0 \leq L^*,$$

so $\theta \in \Theta_{j_n}$ and (44) holds. To bound $\pi_{j_n}(B_{d_{j_n},j_n}(\theta_1, \varepsilon_n))$ in Assumption 5, note that for all $\varepsilon > 0$ and for all $j$,

(45) $\quad \{\theta \in \Theta_j : \|\theta - \theta_1\| \leq \varepsilon\} \subset \{\theta \in \Theta_j : \|\theta - \theta_1\|_\infty \leq \varepsilon\}.$

By (44) and (45) we have

$$\frac{\pi_{j_n}(B_{d_{j_n},j_n}(\theta_1, \varepsilon_n))}{\pi_{j_n}(B_{d_{j_n},j_n}(\beta_n, \varepsilon_n))} \leq \left(\frac{\varepsilon_n}{\varepsilon_n/(m_{j_n}\sqrt{m_{j_n} - 1})}\right)^{m_{j_n}-1} \leq (m_{j_n}\sqrt{m_{j_n} - 1})^{m_{j_n}}.$$

Since
$$\left(\frac{m_{j_n}\sqrt{m_{j_n} - 1}}{A_{j_n}^3}\right)^{m_{j_n}} \leq \left(\frac{m_{j_n}^{1.5}}{(\sqrt{m_{j_n}})^3}\right)^{m_{j_n}} = 1,$$

Assumption 5 holds with $b_2 = 3$ and $K_5 = 1$.

It is clear that Assumption 6 holds with the above $\varepsilon_n$, which tends to zero at the rate $(\log n)^{1/2}n^{-\alpha/(1+2\alpha)}$. By Theorem 1, the result in Lemma 3 holds.

4.4. *Proof of Theorem* 2. We prove Theorem 2 by giving bounds for $U_n$ and $V_n$, and then combining the bounds to show that $U_n/V_n$ converges to zero.

To bound $U_n$, we will use Lemma 9, which is the regression version of Lemma 7.

LEMMA 9. *Suppose that Assumption 7 holds and $\gamma \in (0, 0.25)$ is defined so that*

$$0.0056 = \frac{0.13}{c_{2,c_0,M}\sqrt{c_{1,M,\sigma}}} \frac{\gamma}{\sqrt{1-4\gamma}}.$$

*Then for all $j$ and for all $\xi_j$ such that*

$$\frac{\xi_j}{m_j} \geq \frac{4}{c_{1,M,\sigma}(1-4\gamma)} \log(1072.5 A_j),$$

$$P_{f_o}^* \left[ \frac{1}{n} \sum_{i=1}^n (Y_i - f_o(X_i))^2 - \frac{1}{n} \sum_{i=1}^n (Y_i - f_{\theta,j}(X_i))^2 \right.$$

$$\geq -\gamma \|f_o - f_{\theta,j}\|_{L_2(\mu_X)}^2 + \frac{\xi_j}{n} + 0.0224 \left| \frac{1}{n} \sum_{i=1}^n \varepsilon_i \right| \sqrt{\frac{\xi_j}{n}}$$

$$\left. \text{for some } \theta \in \Theta_j \text{ and } \frac{1}{n} \sum_{i=1}^n |\varepsilon_i| \leq c_0, \frac{1}{n} \sum_{i=1}^n \varepsilon_i^2 \leq c_0^2 \right]$$

$$\leq 15.1 \exp\left(-\frac{c_{1,M,\sigma}(1-4\gamma)\xi_j}{8}\right),$$

*where*

$$c_{1,M,\sigma} = \min\left(\frac{1-\exp(-M^2/(2\sigma^2))}{2M^2}, \frac{1}{2\sigma^2}\right) \quad \text{and} \quad c_{2,c_0,M} = 2(c_0 + 2M).$$

The proof of Lemma 9 is long and is deferred to Section 4.4.1.

Now suppose that Assumption 7 holds. Take $c_0 = 2\sigma$ and define $\gamma$ as in Lemma 9. Let $C_j \geq 0$ be such that $\sum_j e^{-C_j} \leq 1$ and define $\eta_j$ and $a_j$ as (18) and (17), respectively. We will apply Lemma 9 to prove (46), which gives an upper bound for $U_n$,

(46) $$P_{f_o}\left[U_n \leq \alpha \exp\left(\frac{0.0056 Z_n^2}{\sigma} - \frac{\gamma n s_n^2}{4\sigma^2}\right)\right] \geq 1 - (p_1 + p_2 + p_3),$$

where

$$Z_n = \frac{1}{\sqrt{n}\sigma} \sum_{i=1}^n \varepsilon_i \sim N(0,1),$$

$$p_1 = P\left[\frac{1}{n} \sum_{i=1}^n |\varepsilon_i| > c_0\right], \qquad p_2 = P\left[\frac{1}{n} \sum_{i=1}^n \varepsilon_i^2 > c_0^2\right]$$

and

$$p_3 = 15.1 \exp\left(-\frac{c_{1,M,\sigma}(1-4\gamma)\gamma n s_n^2}{32(0.5 + 0.0056\sigma)}\right).$$

To prove (46), take

$$\xi_j = \eta_j + \frac{\gamma n s_n^2}{4(0.5 + 0.0056\sigma)}.$$

Since $U_n$ is

$$\sum_j a_j \int_{(B_{L_2(\mu_X)}, \Theta_j(s_n))^c} \frac{\exp(1/(2\sigma^2) \sum_{i=1}^n (Y_i - f_o(X_i))^2)}{\exp(1/(2\sigma^2) \sum_{i=1}^n (Y_i - f_{\theta,j}(X_i))^2)} \, d\pi_j(\theta),$$

Lemma 9 gives

$$U_n \leq \sum_j a_j \exp\left(\frac{1}{2\sigma^2}\left(-\gamma n s_n^2 + \xi_j + 0.0224 \left|\frac{1}{\sqrt{n}} \sum_{i=1}^n \varepsilon_i\right| \sqrt{\xi_j}\right)\right)$$

$$= \sum_j a_j \exp\left(-\frac{\gamma n s_n^2}{2\sigma^2} + \frac{\xi_j}{2\sigma^2} + \frac{0.0112}{\sigma} |Z_n| \sqrt{\xi_j}\right)$$

$$\leq \sum_j a_j \exp\left(-\frac{\gamma n s_n^2}{2\sigma^2} + \frac{\xi_j}{2\sigma^2} + \frac{0.0056}{\sigma}(Z_n^2 + \xi_j)\right)$$

$$= \exp\left(\frac{0.0056 Z_n^2}{\sigma} - \frac{\gamma n s_n^2}{4\sigma^2}\right) \sum_j a_j \exp\left(\frac{0.5 + 0.0056\sigma}{\sigma^2} \eta_j\right)$$

$$= \alpha \exp\left(\frac{0.0056 Z_n^2}{\sigma} - \frac{\gamma n s_n^2}{4\sigma^2}\right) \sum_j e^{-\eta_j}$$

$$\leq \alpha \exp\left(\frac{0.0056 Z_n^2}{\sigma} - \frac{\gamma n s_n^2}{4\sigma^2}\right)$$

except on a set of probability no greater than

$$P\left[\frac{1}{n}\sum_{i=1}^n |\varepsilon_i| > c_0\right] + P\left[\frac{1}{n}\sum_{i=1}^n \varepsilon_i^2 > c_0^2\right] + 15.1 \sum_j \exp\left(-\frac{c_{1,M,\sigma}(1-4\gamma)\xi_j}{8}\right).$$

Note that

$$\sum_j \exp\left(-\frac{c_{1,M,\sigma}(1-4\gamma)\xi_j}{8}\right)$$

$$= \exp\left(-\frac{c_{1,M,\sigma}(1-4\gamma)\gamma n s_n^2}{32(0.5 + 0.0056\sigma)}\right) \sum_j \exp\left(-\frac{c_{1,M,\sigma}(1-4\gamma)\eta_j}{8}\right)$$

$$\leq \exp\left(-\frac{c_{1,M,\sigma}(1-4\gamma)\gamma n s_n^2}{32(0.5 + 0.0056\sigma)}\right),$$

so now we have the following bound for $U_n$:

$$P_{f_o}\left[U_n \leq \alpha \exp\left(\frac{0.0056 Z_n^2}{\sigma} - \frac{\gamma n s_n^2}{4\sigma^2}\right)\right] \geq 1 - (p_1 + p_2 + p_3).$$

The process of deriving a bound for $V_n$ is the same as in Section 4.1 except for the following changes:

1. Replace $f_o$ by $p_{f_o}$, $f_{\theta,j_n}$ by $p_{f_{\theta,j_n}}$ and Assumptions 2 and 3 by Assumptions 8 and 9.
2. The proof of (28) is modified as follows. First, note that in our regression setting, for all $\theta \in \Theta_j$ and for all $j$,

(47)
$$D(p_{f_o}\|p_{f_{\theta,j}}) = \frac{\|f_o - f_{\theta,j}\|^2_{L_2(\mu_X)}}{2\sigma^2}$$

and

(48)
$$V(p_{f_o}\|p_{f_{\theta,j}}) = \frac{\|f_o - f_{\theta,j}\|^2_{L_2(\mu_X)}}{\sigma^2} + \frac{1}{4\sigma^4}\int (f_o - f_{\theta,j})^4$$
$$\leq \left(\frac{1}{\sigma^2} + \frac{M^2}{\sigma^4}\right)\|f_o - f_{\theta,j}\|^2_{L_2(\mu_X)}.$$

By (47), (48) and (20), for $\theta \in B_{d_{j_n},j_n}(\beta_n, \varepsilon_n)$, we have

$$D(p_{f_o}\|p_{f_{\theta,j_n}}) \leq D(p_{f_o}\|p_{f_{\beta_n,j_n}}) + \frac{\|f_{\beta,j_n} - f_{\theta,j_n}\|^2_{L_2(\mu_X)}}{2\sigma^2}$$
$$\leq \varepsilon_n^2 + \frac{K_0' \varepsilon_n^2}{2\sigma^2}$$

and

$$V(p_{f_o}\|p_{f_{\theta,j_n}}) \leq \left(2 + \frac{2M^2}{\sigma^2}\right) D(p_{f_o}\|p_{f_{\theta,j_n}}).$$

Therefore, (28) holds for

$$t_n^2 = \left(2 + \frac{2M^2}{\sigma^2}\right)\left(1 + \frac{K_0'}{2\sigma^2}\right)\varepsilon_n^2 \stackrel{\text{def}}{=} K'\varepsilon_n^2.$$

3. The process of deriving a lower bound for $V_n$ in (30) is modified as follows:

$$V_n \geq \frac{1}{2} e^{-2nt_n^2} a_{j_n} \pi_{j_n}(B_{D,j_n}(t_n))$$
$$\geq \frac{\alpha e^{-2nt_n^2}}{2} \exp\left(-\left(1 + \frac{1}{2\sigma^2} + \frac{0.0056}{\sigma}\right)\eta_{j_n}\right)\left(\frac{1}{A_{j_n}^{b_1+b_2} K_4 K_5}\right)^{m_{j_n}}$$
$$\geq \frac{\alpha}{2}\exp\left(-2nt_n^2 - \eta_{j_n}\left(1 + \frac{1}{2\sigma^2} + \frac{0.0056}{\sigma}\right.\right.$$

$$
\begin{aligned}
(49) \qquad &\qquad\qquad\qquad + c_1(b_1 + b_2 + (\log(K_4 K_5))_+)\bigg)\bigg) \\
&\stackrel{(19)}{\geq} \frac{\alpha}{2} \exp\bigg(-2nt_n^2 - n\varepsilon_n^2 \bigg(1 + \frac{1}{2\sigma^2} + \frac{0.0056}{\sigma} \\
&\qquad\qquad\qquad + c_1(b_1 + b_2 + (\log(K_4 K_5))_+)\bigg)\bigg) \\
&\geq \frac{\alpha}{2} e^{-Kn\varepsilon_n^2},
\end{aligned}
$$

where $c_1 = c_{1,M,\sigma}$ and

$$K = 2K' + 1 + \frac{1}{2\sigma^2} + \frac{0.0056}{\sigma} + c_1(b_1 + b_2 + (\log(K_4 K_5))_+).$$

Here we have used the fact that

$$\frac{c_1 \eta_j}{m_j} \geq \frac{4}{1 - 4\gamma} \log\left(\frac{1072.5 A_j \sqrt{1 - 4\gamma}}{\gamma}\right) \geq \max(1, \log A_j)$$

for all $j$.

Now we will bound $U_n/V_n$ by combining (46) and (50). In (46), set

$$s_n^2 = \frac{8\sigma^2 K \varepsilon_n^2}{\gamma}.$$

Then

$$
\begin{aligned}
\tilde{\pi}(\tilde{B}_{L_2(\mu_X)}(s_n)^c | X_1, \ldots, X_n) &= \frac{U_n}{V_n} \\
&\leq 2 \exp\left(\frac{0.0056 Z_n^2}{\sigma}\right) \exp(-Kn\varepsilon_n^2)
\end{aligned}
$$

except on a set of probability no greater than

$$p_1 + p_2 + 15.1 \exp\left(-\frac{c_{1,M,\sigma}(1 - 4\gamma) 8\sigma^2 K n \varepsilon_n^2}{32(0.5 + 0.0056\sigma)}\right) + \frac{2}{K' n \varepsilon_n^2},$$

where

$$Z_n = \frac{1}{\sqrt{n}\sigma} \sum_{i=1}^n \varepsilon_i \sim N(0, 1),$$

$$p_1 = P\left[\frac{1}{n} \sum_{i=1}^n |\varepsilon_i| > c_0\right]$$

and

$$p_2 = P\left[\frac{1}{n} \sum_{i=1}^n \varepsilon_i^2 > c_0^2\right].$$

Note that $c_0 = 2\sigma > \max(E|\varepsilon_i|, E\varepsilon_i^2)$, so $p_1 + p_2 \to 0$ as $n \to \infty$. Since $2e^{0.0056 Z_n^2/\sigma}$ converges in distribution and $e^{-Kn\varepsilon_n^2}$ converges to zero by Assumption 6, we have that $2e^{0.0056 Z_n^2/\sigma} e^{-Kn\varepsilon_n^2}$ converges to zero in probability. Therefore, $\tilde{\pi}(\tilde{B}_{L_2(\mu_X)}(s_n)^c | X_1, \ldots, X_n)$ converges to zero in probability as stated in Theorem 2.

4.4.1. *An exponentional inequality.* We claim that to prove Lemma 9, it suffices to prove Lemma 10, which has a slightly different assumption.

ASSUMPTION 10. For some $j \in J$, for $\theta \in \Theta_j$, $\|f_{\theta,j}\|_\infty \leq M$, and there exist constants $A > 0$, $m \geq 1$ and $0 < \rho \leq A$ such that for any $r > 0$, $\delta \leq \rho r$, $\theta \in \Theta_j$, the $\delta$-covering number

$$N(B_{L_2(\mu_X), \Theta_j}(r), \delta, d_{j,\infty}) \leq \left(\frac{Ar}{\delta}\right)^m,$$

where $B_{L_2(\mu_X), \Theta_j}(r) = \{\theta \in \Theta_j : \|f_o - f_{\theta,j}\|_{L_2(\mu_X)} \leq r\}$ and for $\eta, \theta \in \Theta_j$, $d_{j,\infty}(\eta, \theta) = \|f_{\eta,j} - f_{\theta,j}\|_\infty$.

LEMMA 10. *Suppose that Assumption 10 holds with*

$$\rho \geq \frac{0.13}{c_{2,c_0,M}\sqrt{c_{1,M,\sigma}}} \frac{\gamma}{\sqrt{1-4\gamma}}.$$

*Then for $\xi$ such that*

$$\frac{\xi}{m} \geq \frac{4}{c_{1,M,\sigma}(1-4\gamma)} \log\left(15.4 c_{2,c_0,M}\sqrt{c_{1,M,\sigma}} A \frac{\sqrt{1-4\gamma}}{\gamma}\right),$$

$$P^*\left[\frac{1}{n}\sum_{i=1}^n (Y_i - f_o(X_i))^2 - \frac{1}{n}\sum_{i=1}^n (Y_i - f_{\theta,j}(X_i))^2\right.$$

$$\geq -\gamma \|f_o - f_{\theta,j}\|_{L_2(\mu_X)}^2 + \frac{\xi}{n} + 4\left|\frac{1}{n}\sum_{i=1}^n \varepsilon_i\right|\delta$$

$$\left. \text{for some } \theta \in \Theta_j \text{ and } \frac{1}{n}\sum_{i=1}^n |\varepsilon_i| \leq c_0, \frac{1}{n}\sum_{i=1}^n \varepsilon_i^2 \leq c_0^2 \right]$$

$$\leq 15.1 \exp\left(-\frac{c_{1,M,\sigma}(1-4\gamma)\xi}{8}\right),$$

*where*

$$\delta = \frac{2\gamma}{15.4 c_{2,c_0,M}\sqrt{c_{1,M,\sigma}}(1-4\gamma)} \sqrt{\frac{\xi}{n}},$$

$$c_{1,M,\sigma} = \min\left(\frac{1-\exp(-M^2/(2\sigma^2))}{2M^2}, \frac{1}{2\sigma^2}\right) \quad and \quad c_{2,c_0,M} = 2(c_0 + 2M).$$

To see that the claim is true, note that in the proof for (26), $d_H$ can be replaced by $L_2(\mu_X)$. Therefore, if Assumption 7 holds, then for all $j \in J$, Assumption 10 holds with $A = 3A_j$ and $\rho = 0.0056$. Suppose that Lemma 10 is true. Then Lemma 9 follows by setting $\rho = 0.0056$ and choosing $\gamma$ such that

$$\rho = \frac{0.13}{c_{2,c_0,M}\sqrt{c_{1,M,\sigma}}} \frac{\gamma}{\sqrt{1-4\gamma}}.$$

PROOF OF LEMMA 10. We follow the proof of Lemma 0 in Yang and Barron (1998). First, divide the space $\Theta_j$ into rings

$$\Theta_{j,i} = \{\theta \in \Theta_j : r_{i-1} \leq \|f_o - f_{\theta,j}\|_{L_2(\mu_X)} \leq r_i\}, \qquad i = 0, 1, \ldots,$$

where $r_i = 2^{i/2}\sqrt{\xi/n}$ for $i \geq 0$ and $r_{-1} = 0$. For each ring $\Theta_{j,i}$, we will use a chaining argument to bound

$$q_i \stackrel{\text{def}}{=} P^*\left[\frac{1}{n}\sum_{i'=1}^n (Y_{i'} - f_o(X_{i'}))^2 - \frac{1}{n}\sum_{i'=1}^n (Y_{i'} - f_{\theta,j}(X_{i'}))^2\right.$$

$$\geq -\gamma\|f_o - f_{\theta,j}\|_{L_2(\mu_X)}^2 + \frac{\xi}{n} + 4\left|\frac{1}{n}\sum_{i'=1}^n \varepsilon_{i'}\right|\delta$$

$$\left. \text{for some } \theta \in \Theta_{j,i} \text{ and } \frac{1}{n}\sum_{i'=1}^n |\varepsilon_{i'}| \leq c_0, \frac{1}{n}\sum_{i'=1}^n \varepsilon_{i'}^2 \leq c_0^2\right].$$

Then we will put all the bounds for $q_i$ together to complete the proof. So let us focus on one $\Theta_{j,i}$ first. Let $\{\delta_k\}_{k=0}^\infty$ be a sequence decreasing to zero with $\delta_0 \leq \min(\rho r_0, \delta)$ and define $\tilde{\delta}_k = \delta_k$ for $k \geq 1$ and $\tilde{\delta}_0 = \delta_0/2$. Then by assumption we can find a sequence of nets $\tilde{F}_0, \tilde{F}_1, \ldots$, where each $\tilde{F}_k$ is a $\tilde{\delta}_k$ net in $\Theta_{j,i}$ satisfying the cardinal number constraint in Assumption 10. In other words, for each $k$, there exists a mapping $\tilde{\tau}_k : \Theta_{j,i} \to \tilde{F}_k$ such that $\|f_{\tilde{\tau}_k(\theta),j} - f_{\theta,j}\|_\infty \leq \tilde{\delta}_k$ for all $\theta \in \Theta_{j,i}$, and

$$\text{card}(\tilde{F}_k) \leq \left(\frac{Ar_k}{\tilde{\delta}_k}\right)^m.$$

Instead of applying the chaining argument using the nets $\tilde{F}_k$, we will modify the net $\tilde{F}_0$ first and then apply the chaining argument using the nets $F_k$, where $F_k = \tilde{F}_k$ for $k \geq 1$ and $F_0$ is the modified $\tilde{F}_0$. Now modify the net $\tilde{F}_0$ in the following way: Consider a positive number $\varepsilon$. For each $\tilde{\theta}_0$ in $\tilde{F}_0$, find $\theta_0$ in

$$\tilde{\tau}_0^{-1}(\tilde{\theta}_0) = \{\theta \in \Theta_{j,i} : \tilde{\tau}_0(\theta) = \tilde{\theta}_0\}$$

such that
$$\|f_o - f_{\tilde\theta_0,j}\|^2_{L_2(\mu_X)} \le \inf_{\theta\in\tilde\tau_0^{-1}(\tilde\theta_0)} \|f_o - f_{\theta,j}\|^2_{L_2(\mu_X)} + \varepsilon.$$

Define $\tau(\tilde\theta_0) = \theta_0$, and $F_0 = \{\tau(\tilde\theta_0) : \tilde\theta_0 \in \tilde F_0\}$. Define $\tau_0 = \tau(\tilde\tau_0)$ and $\tau_k = \tilde\tau_k$ for $k \ge 1$. Then by the triangle inequality, $\|f_{\tau_0(\theta),j} - f_{\theta,j}\|_\infty \le \delta_0$, so $F_0$ is a $\delta_0$ net and for each $k$, $F_k$ is a $\delta_k$ net. Now we can start the chaining argument. For each $\theta \in \Theta_{j,i}$, define

$$l_0 = \frac{1}{n}\sum_{i=1}^n (Y_i - f_o(X_i))^2 - \frac{1}{n}\sum_{i=1}^n (Y_i - f_{\tau_0(\theta),j}(X_i))^2$$

and

$$l_k = \frac{1}{n}\sum_{i=1}^n (Y_i - f_{\tau_{k-1}(\theta),j}(X_i))^2 - \frac{1}{n}\sum_{i=1}^n (Y_i - f_{\tau_k(\theta),j}(X_i))^2$$

for $k \ge 1$. Then

$$\frac{1}{n}\sum_{i=1}^n (Y_i - f_o(X_i))^2 - \frac{1}{n}\sum_{i=1}^n (Y_i - f_{\theta,j}(X_i))^2 = l_0 + \sum_{k=1}^\infty l_k.$$

Now, instead of giving bounds for $l_k - El_k$ as in Yang and Barron (1998), we will give bounds for $l_k - E_\varepsilon l_k$, where

$$E_\varepsilon l_k = \frac{2}{n}\sum_{i=1}^n \varepsilon_i \int (f_{\tau_k(\theta),j} - f_{\tau_{k-1}(\theta),j})\, d\mu_X$$
$$+ \|f_o - f_{\tau_{k-1}(\theta),j}\|^2_{L_2(\mu_X)} - \|f_o - f_{\tau_k(\theta),j}\|^2_{L_2(\mu_X)}$$

is the conditional expectation of $l_k$ given $\varepsilon_1, \ldots, \varepsilon_n$ for $k \ge 1$. Note that

$$\sum_{k=1}^\infty E_\varepsilon l_k = 2\left(\frac{1}{n}\sum_{i=1}^n \varepsilon_i\right)\int (f_{\theta,j} - f_{\tau_0(\theta),j})\mu_X$$
$$+ \|f_o - f_{\tau_0(\theta),j}\|^2_{L_2(\mu_X)} - \|f_o - f_{\theta,j}\|^2_{L_2(\mu_X)}$$
$$\le 2\left|\frac{1}{n}\sum_{i=1}^n \varepsilon_i\right|\left|\int (f_{\theta,j} - f_{\tau_0(\theta),j})\mu_X\right| + \varepsilon$$
$$\le 4\left|\frac{1}{n}\sum_{i=1}^n \varepsilon_i\right|\delta_0 + \varepsilon \le 4\left|\frac{1}{n}\sum_{i=1}^n \varepsilon_i\right|\delta + \varepsilon,$$

so

$$q_i \le P^*(B_0 \cap B)$$
$$\le P^*\left(\left\{l_0 \ge -2\gamma r_i^2 + \frac{\xi}{n} - \varepsilon \text{ for some } \theta \in \Theta_{j,i}\right\} \cap B\right)$$

$$+ \sum_{k=1}^{\infty} P^*(\{l_k - E_\varepsilon l_k \geq \eta_k \text{ for some } \theta \in \Theta_{j,i}\} \cap B)$$

$$\stackrel{\text{def}}{=} q_i^{(1)} + \sum_{k=1}^{\infty} q_{i,k}^{(2)}$$

if $\sum_{k=1}^{\infty} \eta_k \leq \gamma r_i^2$, where

$$B_0 = \left\{ l_0 + \sum_{k=1}^{\infty}(l_k - E_\varepsilon l_k) \geq -\varepsilon - \gamma\|f_o - f_{\theta,j}\|_{L_2(\mu_X)}^2 + \frac{\xi}{n} \text{ for some } \theta \in \Theta_{j,i} \right\}$$

and

$$B = \left\{ \frac{1}{n}\sum_{i=1}^{n} |\varepsilon_i| \leq c_0, \frac{1}{n}\sum_{i=1}^{n} \varepsilon_i^2 \leq c_0^2 \right\}. \qquad \square$$

To bound $q_i^{(1)}$, we will use the following inequality of Chernoff (1952):

FACT 3. Suppose that $X_i$ are i.i.d. from a distribution with density $g_2$ with respect to measure $\mu$ and $g_1$ is a density with respect to the same measure. Then

$$P\left[\frac{1}{n}\sum_{i=1}^{n} \log \frac{g_1(X_i)}{g_2(X_i)} \geq t\right] \leq \exp\left(-\frac{n}{2}(d_H^2(g_1, g_2) + t)\right).$$

Since

$$l_0 = \frac{2\sigma^2}{n} \sum_{i=1}^{n} \log \frac{p_{f_{\tau_0(\theta),j}}(X_i)}{p_{f_o}(X_i)},$$

Fact 3 implies that for a $\tau_0(\theta)$,

$$(50) \qquad P[l_0 \geq t] \leq \exp\left(-\frac{n}{2}(d_H^2(p_{f_{\tau_0(\theta),j}}, p_{f_o}) + t/(2\sigma^2))\right).$$

To replace the Hellinger distance $d_H^2(p_{f_{\tau_0(\theta),j}}, p_{f_o})$ with the $L^2$ distance $\|f_{\tau_0(\theta),j} - f_o\|_{L_2(\mu_X)}$ in (50), note that

$$d_H^2(p_{f_{\tau_0(\theta),j}}, p_{f_o}) = 2\int \left(1 - \exp\left(-\frac{(f_{\tau_0(\theta),j}(x) - f_o(x))^2}{8\sigma^2}\right)\right) d\mu(x)$$

$$(51) \qquad \geq \frac{1 - \exp(-M^2/(2\sigma^2))}{2M^2} \int (f_{\tau_0(\theta),j}(x) - f_o(x))^2 d\mu(x)$$

$$\stackrel{\text{def}}{=} c_{0,M,\sigma} \|f_{\tau_0(\theta),j} - f_o\|_{L_2(\mu_X)}^2.$$

Here the equality follows from direct calculation and the inequality follows from the fact that $(1 - e^{-x})/x$ is decreasing with $x$ on $(0, \infty)$ and that $\|f_{\tau_0(\theta),j}\|_\infty, \|f_o\|_\infty \leq M$. Now by (50) and (51), we have

$$P[l_0 \geq t] \leq \exp\left(-\frac{n}{2}(c_{0,M,\sigma}\|f_{\tau_0(\theta),j} - f_o\|^2_{L_2(\mu_X)} + t/(2\sigma^2))\right)$$

$$\leq \exp\left(-\frac{c_{1,M,\sigma}n}{2}(\|f_{\tau_0(\theta),j} - f_o\|^2_{L_2(\mu_X)} + t)\right),$$

where $c_{1,M,\sigma} = \min(c_{0,M,\sigma}, 1/(2\sigma^2))$. Set $t = -2\gamma r_i^2 + \frac{\xi}{n} - \varepsilon$. Then for a $\tau_0(\theta)$,

$$P\left[l_0 \geq -2\gamma r_i^2 + \frac{\xi}{n} - \varepsilon\right] \leq \exp\left(-\frac{c_{1,M,\sigma}n}{2}\left(r_{i-1}^2 - 2\gamma r_i^2 + \frac{\xi}{n} - \varepsilon\right)\right).$$

Therefore,

$$\begin{aligned}
q_i^{(1)} &\leq \text{card}(F_0)\exp\left(-\frac{c_{1,M,\sigma}n}{2}\left(r_{i-1}^2 - 2\gamma r_i^2 + \frac{\xi}{n} - \varepsilon\right)\right) \\
&\leq \text{card}(F_0)\exp\left(-\frac{c_{1,M,\sigma}n}{2}\left((i+1)(1-4\gamma)\frac{\xi}{n} - \varepsilon\right)\right),
\end{aligned} \tag{52}$$

where the last inequality was verified in Yang and Barron (1998), from the end of page 111 to the beginning of page 112.

To bound $q_{i,k}^{(2)}$, we will use Hoeffding's inequality.

FACT 4. *Suppose that $\{Y_i\}_{i=1}^n$ are independent with mean zero and that $a_i \leq Y_i \leq b_i$ for all $i$. Then for $\eta > 0$,*

$$P\left[\sum_{i=1}^n Y_i \geq \eta\right] \leq \exp\left(\frac{-2\eta^2}{\sum_{i=1}^n (b_i - a_i)^2}\right).$$

For a pair $(\tau_{k-1}(\theta), \tau_k(\theta))$,

$$\begin{aligned}
|(Y_i &- f_{\tau_{k-1}(\theta),j}(X_i))^2 - (Y_i - f_{\tau_k(\theta),j}(X_i))^2| \\
&\leq 2|f_{\tau_{k-1}(\theta),j}(X_i) - f_{\tau_k(\theta),j}(X_i)| \\
&\quad \times \left|\varepsilon_i + f_o(X_i) - \frac{f_{\tau_{k-1}(\theta),j}(X_i) + f_{\tau_k(\theta),j}(X_i)}{2}\right| \\
&\leq 2(\delta_{k-1} + \delta_k)(|\varepsilon_i| + 2M) \leq 4(|\varepsilon_i| + 2M)\delta_{k-1}.
\end{aligned}$$

By Hoeffding's inequality, the conditional probability

$$P[l_k - E_\varepsilon l_k \geq \eta | \varepsilon_1, \ldots, \varepsilon_n] \leq \exp\left(\frac{-2n^2\eta^2}{\sum_{i=1}^n 64(|\varepsilon_i| + 2M)^2 \delta_{k-1}^2}\right)$$

$$\leq \exp\left(\frac{-2n\eta^2}{64(c_0 + 2M)^2 \delta_{k-1}^2}\right)$$

if $\sum_{i=1}^{n} |\varepsilon_i|/n \leq c_0$ and $\sum_{i=1}^{n} \varepsilon_i^2/n \leq c_0^2$. Integrating the conditional probability over set $B$, we have

$$P(\{l_k - E_\varepsilon l_k \geq \eta\} \cap B) \leq \exp\left(\frac{-2n\eta^2}{64(c_0 + 2M)^2 \delta_{k-1}^2}\right).$$

Therefore,

(53) $$q_{i,k}^{(2)} \leq \operatorname{card}(F_{k-1}) \operatorname{card}(F_k) \exp\left(\frac{-2n\eta_k^2}{64(c_0 + 2M)^2 \delta_{k-1}^2}\right).$$

Now combine (52) and (53) and let $\varepsilon \to 0$. Then we have

$$q_i \leq \operatorname{card}(F_0) \exp\left(-\frac{nc_{1,M,\sigma}}{2}(i+1)(1-4\gamma)\frac{\xi}{n}\right)$$

$$+ \sum_{k=1}^{\infty} \operatorname{card}(F_{k-1}) \operatorname{card}(F_k) \exp\left(\frac{-2n\eta_k^2}{64(c_0 + 2M)^2 \delta_{k-1}^2}\right)$$

$$\leq \left(\frac{Ar_i}{\tilde{\delta}_0}\right)^m \exp\left(-\frac{c_{1,M,\sigma}}{2}(i+1)(1-4\gamma)\xi\right)$$

$$+ \sum_{k=1}^{\infty} \left(\frac{Ar_i}{\tilde{\delta}_{k-1}}\right)^m \left(\frac{Ar_i}{\tilde{\delta}_k}\right)^m \exp\left(\frac{-2n\eta_k^2}{64(c_0 + 2M)^2 \delta_{k-1}^2}\right).$$

Now choose $\delta_0, \delta_k$ so that

$$\log\left(\frac{Ar_0}{\tilde{\delta}_k}\right)^m = \frac{c_{1,M,\sigma}(k+1)(1-4\gamma)\xi}{4}$$

and $\eta_k$ such that

$$\frac{2n\eta_k^2}{64(c_0 + 2M)^2 \delta_{k-1}^2}$$
$$= im\log 2 + \frac{(2k+1)c_{1,M,\sigma}(1-4\gamma)\xi}{4} + \frac{(i+1)kc_{1,M,\sigma}(1-4\gamma)\xi}{8}.$$

Now the bound for $q_i$ becomes

$$q_i \leq 2^{im/2} \exp\left(\frac{c_{1,M,\sigma}(1-4\gamma)\xi}{4}\right) \exp\left(-\frac{c_{1,M,\sigma}}{2}(i+1)(1-4\gamma)\xi\right)$$

$$+ \sum_{k=1}^{\infty} \exp\left(-\frac{(i+1)c_{1,M,\sigma}k(1-4\gamma)\xi}{8}\right)$$

$$\leq \exp\left(\frac{im}{2}\log 2 - \frac{(i+1)c_{1,M,\sigma}(1-4\gamma)\xi}{4}\right)$$

$$+ \exp\left(-\frac{(i+1)c_{1,M,\sigma}(1-4\gamma)\xi}{8}\right)$$

$$\times \left(1 - \exp\left(-\frac{(i+1)c_{1,M,\sigma}(1-4\gamma)\xi}{8}\right)\right)^{-1}.$$

Note that by assumption,

$$\frac{m}{2}\log\frac{2A}{\rho_0} \leq \frac{c_{1,M,\sigma}(1-4\gamma)\xi}{8},$$

where

$$\rho_0 = \frac{2\gamma}{15.4 c_{2,c_0,M}\sqrt{c_{1,M,\sigma}(1-4\gamma)}}.$$

Since $\rho_0 \leq \rho \leq A$, we have

$$(54) \quad \frac{\log 2}{2} \leq \frac{m}{2}\log 2 \leq \frac{m}{2}\log\frac{2A}{\rho_0} \leq \frac{c_{1,M,\sigma}(1-4\gamma)\xi}{8},$$

so

$$q_i \leq \exp\left(-\frac{(i+1)c_{1,M,\sigma}(1-4\gamma)\xi}{8}\right)$$
$$\times \left(1 + \left(1 - \exp\left(-\frac{c_{1,M,\sigma}(1-4\gamma)\xi}{8}\right)\right)^{-1}\right)$$
$$\leq \left(1 + \frac{\sqrt{2}}{\sqrt{2}-1}\right)\exp\left(-\frac{(i+1)c_{1,M,\sigma}(1-4\gamma)\xi}{8}\right)$$

and

$$P^*\left[\frac{1}{n}\sum_{i=1}^n (Y_i - f_o(X_i))^2 - \frac{1}{n}\sum_{i=1}^n (Y_i - f_{\theta,j}(X_i))^2\right.$$
$$\geq -\gamma\|f_o - f_{\theta,j}\|^2_{L_2(\mu_X)} + \frac{\xi}{n} + 4\left|\frac{1}{n}\sum_{i=1}^n \varepsilon_i\right|\delta$$
$$\left. \text{for some } \theta \in \Theta_j \text{ and } \frac{1}{n}\sum_{i=1}^n |\varepsilon_i| \leq c_0, \frac{1}{n}\sum_{i=1}^n \varepsilon_i^2 \leq c_0^2 \right]$$
$$\leq \sum_{i=0}^\infty q_i$$
$$\leq \left(1 + \frac{\sqrt{2}}{\sqrt{2}-1}\right)\exp\left(-\frac{c_{1,M,\sigma}(1-4\gamma)\xi}{8}\right)$$
$$\times \left(1 - \exp\left(-\frac{c_{1,M,\sigma}(1-4\gamma)\xi}{8}\right)\right)^{-1}$$
$$\stackrel{(54)}{\leq} 15.1\exp\left(-\frac{c_{1,M,\sigma}(1-4\gamma)\xi}{8}\right).$$

It remains to check that $\{\delta_k\}_{k=0}^{\infty}$ is a decreasing sequence

$$\sum_{k=1}^{\infty} \eta_k \leq \gamma r_i^2, \tag{55}$$

and

$$\delta_0 \leq \min(r_0 \rho, \delta), \tag{56}$$

as claimed in the beginning of the proof. By (54), $\delta_0/\delta_1 \geq 1$, so $\{\delta_k\}_{k=0}^{\infty}$ is decreasing by construction. To verify (55), let $c_2 = 2(c_0 + 2M)$ and $c_1 = c_{1,M,\sigma}$. Then

$$\eta_1 = 2c_2 A \sqrt{\frac{\xi}{n}} \exp\left(-\frac{c_1(1-4\gamma)\xi}{4m}\right) \sqrt{\frac{im8\log 2}{n} + \frac{(i+7)c_1(1-4\gamma)\xi}{n}}$$

$$\stackrel{(54)}{\leq} 2c_2 A \frac{\xi}{n} \sqrt{c_1(1-4\gamma)} \exp\left(-\frac{c_1(1-4\gamma)\xi}{4m}\right) \sqrt{3i+9},$$

and for $k \geq 2$,

$$\eta_k = c_2 A \sqrt{\frac{\xi}{n}} \exp\left(-\frac{c_1 k(1-4\gamma)\xi}{4m}\right)$$

$$\times \sqrt{\frac{im8\log 2}{n} + \frac{2(2k+1)c_1(1-4\gamma)\xi}{n} + \frac{(i+1)kc_1(1-4\gamma)\xi}{n}}$$

$$\stackrel{(54)}{\leq} c_2 A \frac{\xi}{n} \sqrt{c_1(1-4\gamma)} \exp\left(-\frac{c_1 k(1-4\gamma)\xi}{4m}\right)$$

$$\times \sqrt{2(2k+1) + (i+1)(k+2)}$$

$$\leq c_2 A \frac{\xi}{n} \sqrt{c_1(1-4\gamma)} \exp\left(-\frac{c_1 k(1-4\gamma)\xi}{4m}\right) \sqrt{(i+5)(k+2)}$$

$$\leq c_2 A \frac{\xi}{n} \sqrt{c_1(1-4\gamma)} \exp\left(-\frac{c_1 k(1-4\gamma)\xi}{8m}\right) \sqrt{i+5}.$$

Therefore,

$$\sum_{k=1}^{\infty} \eta_k \leq c_2 A \frac{\xi}{n} \sqrt{c_1(1-4\gamma)} \exp\left(-\frac{c_1(1-4\gamma)\xi}{4m}\right) \sqrt{i+5}$$

$$\times \left(2\sqrt{3} + \frac{1}{1 - \exp(-c_1(1-4\gamma)\xi/(8m))}\right)$$

$$\leq c_2 A \frac{\xi}{n} \sqrt{c_1(1-4\gamma)} \exp\left(-\frac{c_1(1-4\gamma)\xi}{4m}\right) \sqrt{5}\, 2^i$$

$$\times \left(2\sqrt{3} + \frac{1}{1 - \exp(-c_1(1-4\gamma)\xi/(8m))}\right)$$

$$\stackrel{(54)}{\leq} c_2 A \frac{\xi}{n} \sqrt{c_1(1-4\gamma)} \exp\left(-\frac{c_1(1-4\gamma)\xi}{4m}\right) \sqrt{5}\, 2^i \left(2\sqrt{3} + \frac{\sqrt{2}}{\sqrt{2}-1}\right)$$

$$\leq 15.4 c_2 \sqrt{c_1} A \frac{\sqrt{1-4\gamma}}{\gamma} \exp\left(-\frac{c_1(1-4\gamma)\xi}{4m}\right) \gamma 2^i \frac{\xi}{n}$$

$$= 15.4 c_2 \sqrt{c_1} A \frac{\sqrt{1-4\gamma}}{\gamma} \exp\left(-\frac{c_1(1-4\gamma)\xi}{4m}\right) \gamma r_i^2.$$

To make (55) hold, it is sufficient to require that

$$\frac{\xi}{m} \geq \frac{4}{c_1(1-4\gamma)} \log\left(15.4 c_2 \sqrt{c_1} A \frac{\sqrt{1-4\gamma}}{\gamma}\right)$$

as in the assumption. Now it remains to verify (56). (56) follows from the fact that

$$\delta_0 = 2A \sqrt{\frac{\xi}{n}} \exp\left(-\frac{c_1(1-4\gamma)\xi}{4m}\right) \stackrel{(54)}{\leq} \rho_0 \sqrt{\frac{\xi}{n}} = \delta$$

and that

$$\frac{\delta_0}{r_0} = 2A \exp\left(-\frac{c_1(1-4\gamma)\xi}{4m}\right) \stackrel{(54)}{\leq} 2A \frac{\rho_0}{2A} \leq \rho.$$

The proof for Lemma 10 is complete. □

4.5. *Proof of Lemma 5.* We will prove Lemma 5 by verifying the assumptions in Theorem 2. To verify Assumption 7, we will apply Lemma 6. Following the same arguments in the verification of Assumption 1 of Lemma 1 in Section 4.2, we have that (8) and (9) hold with $T_1 = 1$ and $T_2 = 1/(\sqrt{q}(2q+1)9^{q-1})$. By Lemma 6, Assumption 7 holds for $A_j$ and $m_j$ in (22). Note that for the $C_j$ specified in (22), $\sum_j e^{-C_j} = e^{-2}/(1-e^{-1})^3 < 1$ as required.

To verify Assumption 8, we choose $j_n$ and $\beta_n$ as in the verification for Assumption 2 in the proof of Lemma 1 except for the following changes:

1. Fact 1 is replaced by Fact 5.

FACT 5. For $j$ such that $q \geq s+1$, there exists $\beta \in R^{m_j}$ such that

(57)
$$\|D^r(f_o - f_{\beta,j})\|_\infty \leq \alpha_q \left(\frac{1}{k+1}\right)^{s-r} M_0 \quad \text{for } 0 \leq r \leq s-1,$$
$$\|D^s f_{\beta,j}\|_\infty \leq \alpha_q M_0,$$

where $M_0 = \max_{0 \leq r \leq s} \|D^r f_o\|_{L_\infty}$.

The above fact follows from (6.50) in Schumaker (1981).

2. $\beta_n \in R^{m_{j_n}}$ is chosen so that

(58) $$\|f_o - f_{\beta_n,j_n}\|_\infty \leq \alpha_{q^*} M_0 \left(\frac{1}{k_n+1}\right)^s.$$

By (47), (48) and (58), for the above $j_n$ and $\beta_n$,

$$\max(D(f_o\|f_{\beta_n,j_n}), V(f_o\|f_{\beta_n,j_n})) + \frac{\eta_{j_n}}{n} \leq c_1 n^{-2s/(1+2s)},$$

so Assumption (2) holds if $\beta_n \in \Theta_{j_n}$ and

(59) $$\varepsilon_n^2 = c_1 n^{-2s/(1+2s)}.$$

To verify that $\beta_n \in \Theta_{j_n}$, we need to make sure $\max_{0 \leq r \leq q-1} \|D^r f_{\beta_n,j_n}\|_{L_\infty} \leq L^*$ and $\|f_{\beta_n,j_n}\|_\infty \leq M$. The first condition follows from the second equation in (57). The second condition holds for large $n$ because of (58) and the fact that $\|f_o\| < M$. Therefore, Assumption 8 holds for large $n$ for the $\varepsilon_n$ in (59).

Assumption 9 holds with $d_{j_n}(\eta,\theta) = \|f_{\eta,j_n} - f_{\theta,j_n}\|_\infty$ for all $\eta, \theta \in \Theta_{j_n}$ since (20) holds with $K_0' = 1$.

For Assumption 4, the verification is the same as the one for Assumption 4 in the proof of Lemma 1.

To verify Assumption 5, we need to bound $\pi_{j_n}(B_{d_{j_n},j_n}(\beta_n, \varepsilon_n))$ by showing that

(60) $$\left\{\theta \in R^{m_{j_n}} : \|\theta - \beta_n\|_\infty \leq c_6 \left(\frac{1}{k_n+1}\right)^s\right\} \subset B_{d_{j_n},j_n}(\beta_n, \varepsilon_n),$$

where $c_6 = \min(1, \sqrt{c_1}/(\sup_n n^{s/(1+2s)}(k_n+1)^{-s}))$. For $\theta \in R^{m_{j_n}}$ such that $\|\theta - \beta_n\|_\infty \leq c_6(1/(k_n+1))^s$, we will prove (37) and (38). The inequality (37) follows from the same arguments as in the verification for (37) in the proof of Lemma 1, except that $\|\log f_{\theta,j_n} - \log f_{\beta_n,j_n}\|_\infty$ is replaced by $\|f_{\theta,j_n} - f_{\beta_n,j_n}\|_\infty$ and the factor 2 is dropped. To prove (38), note that for $0 \leq r \leq s$,

$$\|D^r f_{\theta,j_n}\|_\infty \leq L^* \quad \text{and} \quad \|D^s f_{\theta,j_n}\|_{L_\infty} \leq L^*,$$

where the results follow from the same arguments for the verification of (38) in the proof of Lemma 1 except that $\log f_{\theta,j_n}$ is replaced by $f_{\theta,j_n}$, $\log f_o$ is replaced by $f_o$ and the case $r=0$ is combined with the case $0 < r < s$ here. Also,

$$\begin{aligned}
\|f_{\theta,j_n}\|_\infty &= \|\theta' B\|_\infty \\
&\leq \|\theta' B - \beta_n' B\|_\infty + \|\beta_n' B - f_o\|_\infty + \|f_o\|_\infty \\
&\stackrel{(39),(57)}{\leq} \|\theta - \beta_n\|_\infty + \alpha_{q^*}\left(\frac{1}{k_n+1}\right)^s M_0 + \|f_o\|_\infty \\
&\leq \left(\frac{1}{k_n+1}\right)^s (1 + \alpha_{q^*} M_0) + \|f_o\|_\infty < M
\end{aligned}$$

for large $n$ since $\|f_o\|_\infty < M$. Therefore, $\theta \in \Theta_{k_n,q^*,L^*}$ and (60) holds.

To bound $\pi_{j_n}(B_{d_{j_n},j_n}(\theta_1, \varepsilon_n))$ in Assumption 5, note that by Lemma 4.3 of Ghosal, Ghosh and van der Vaart (2000), there exists $\beta_{q^*}^* > 1$ such that for all $\varepsilon > 0$ and for all $j$,

(61) $\quad \{\theta \in \Theta_j : \|f_{\theta,j} - f_{\theta_1,j}\|_\infty \leq \varepsilon\} \subset \{\theta \in \Theta_j : \|\theta - \theta_1\|_\infty \leq \beta_{q^*}^* \varepsilon\}.$

Then by (61) and (60), following the arguments after the verification of (40) in the proof of Lemma 1, Assumption 5 holds with $K_5 = \beta_{q^*}^*(1 + (c_4\sqrt{c_1})^{1/s})^s/c_6$ and $b_2 = 0$.

For Assumption 6, it should be clear that it holds with the $\varepsilon_n$ specified in (59). Apply Theorem 2 and we have the result in Lemma 5.

**Acknowledgments.** The author thanks Professor L. Wasserman for supervising the author's thesis work, from which this paper is taken. The author also thanks the referees and an Associate Editor for constructive suggestions to improve the readability and accuracy of this paper.

Department of Statistics
Iowa State University
314 Snedecor Hall
Ames, Iowa 50011-1210
USA
e-mail: tmhuang@iastate.edu